\documentclass[11pt, a4paper]{amsart}

\usepackage{latexsym}

\usepackage{color}

\usepackage{amssymb}
\usepackage[centertags]{amsmath}
\usepackage{verbatim}

\usepackage{units}

\usepackage{datetime}

\usepackage{amsthm}

\newtheorem{theorem}{Theorem}

\newtheorem{lemma}{Lemma}
\newtheorem{question}[theorem]{Question}

\newtheorem{corollary}[theorem]{Corollary}

\newtheorem{definition}[equation]{Definition}

\newtheorem{claim}[equation]{Claim}

\theoremstyle{definition}

\newcommand{\theoremname}{testing}
\numberwithin{equation}{section}
\numberwithin{theorem}{section}
\numberwithin{lemma}{section}

\renewcommand{\phi}{\varphi}

\newcommand{\e}{\varepsilon}
\newcommand{\Om}{\Omega}
\newcommand{\om}{\omega}

\newcommand{\la}{\lambda}

\newcommand{\cE}{\mathcal E}

\newcommand{\cH}{\mathcal H}

\newcommand{\cP}{\mathcal P}

\newcommand{\frM}{\mathfrak M}

\def\done{{1\hskip-2.5pt{\rm l}}}
\def\eqdef{\overset{\rm{def}}{=}}
\def\re#1{\mathrm{Re}\,#1}
\def\RFfuncs{L^2_{\tt RF}}

\renewcommand{\le}{\leqslant}
\renewcommand{\ge}{\geqslant}

\newcommand{\bR}{\mathbb R}
\newcommand{\bC}{\mathbb C}
\newcommand{\bZ}{\mathbb Z}
\newcommand{\bD}{\mathbb D}
\newcommand{\bT}{\mathbb T}
\newcommand{\bN}{\mathbb N}
\newcommand{\bQ}{\mathbb Q}

\title[Log-integrability of Rademacher Fourier series]
{Log-integrability of Rademacher Fourier series, \\
with applications to random analytic functions}

\author{Fedor Nazarov}
\address{Department of Mathematical Sciences, Kent State University, Kent OH 44242, USA
}
\email{nazarov@math.kent.edu}

\author{Alon Nishry}
\address{School of Mathematical Sciences\\
Tel Aviv University\\
Tel Aviv 69978\\
Israel}
\email{alonnish@post.tau.ac.il}
\thanks{}

\author{Mikhail Sodin}
\address{School of Mathematical Sciences\\
Tel Aviv University\\
Tel Aviv 69978\\
Israel}
\email{sodin@post.tau.ac.il}
\thanks{This work was partially supported by grant No. 2006136 of the United States - Israel
Binational Science Foundation (F.N., A.N., M.S.), by U.S. National Science Foundation Grant
DMS-0800243 (F.N.), and by grant No. 166/11  of the Israel Science Foundation of the
Israel Academy
of Sciences and Humanities (A.N., M.S.)}

\begin{document}

\begin{abstract}
We prove that any power of the logarithm of Fourier series with random signs
is integrable. This result has applications to the distribution of values
of random Taylor series, one of which answers a
long-standing question by J.-P. Kahane.
\end{abstract}

\maketitle

\date{03.01.2013}

\section{Introduction}

In this work, we consider Rademacher Fourier series
\[
f(\theta) = \sum_{k\in \bZ} \xi_k a_k  e^{2\pi {\rm i} k \theta}
\]
where $\xi_k$ are independent Rademacher random variables, which take the values
$\pm 1$ with probability $\frac{1}{2}$ each,
and random Taylor series
\[
F(z) = \sum_{k\in \bZ_+} \zeta_k  z^k
\]
with independent symmetric complex-valued random variables $\zeta_k$.
Recall that the complex-valued random variable $\zeta$ is called {\em symmetric}
if $-\zeta$ has the same distribution as $\zeta$.
In the Fourier case, the sequence of deterministic complex coefficients
$\bigl\{ a_k \bigr\}$ belongs to $\ell^2(\bZ)$; in the Taylor case,
we assume that the radius of convergence is almost surely (a.s., for short) positive.

\subsection{Some motivation}\label{subsect:intro}
There are several long-standing questions pertaining to the distribution of
values of random Taylor series. For these questions, the Rademacher case already presents main
difficulties. Moreover, in many instances, due to Kahane's
``reduction principle''~\cite[Section~1.7]{Kahane},
the case of more general random symmetric coefficients can be reduced to the
Rademacher case. Here, we explain the central  r\^{o}le played by the
logarithmic integrability of the Rademacher Fourier series in our approach to
some of these questions.

Put $\bT=\bR/\bZ$, and denote by $m$ the normalized Lebesgue measure on $\bT$.
Consider a random Taylor series
\[
F(z) = \sum_{k\ge 0} \xi_k a_k z^k
\]
with independent identically distributed complex-valued random coefficients $\xi_k$
normalized by $\cE |\xi|^2=1$.
Let $R$, $0<R\le\infty$ denote the radius of convergence of this Taylor series.
Note that
\[
\mathcal E\bigl\{ |F(z)|^2 \bigr\} =
\sum_{k \ge 0} |a_k|^2 r^{2k}
\]
for all $z$ with $|z|=r$. We denote the RHS by $\sigma^2_F(r)$.
We will always assume that $\sigma_F(r)\to \infty$ as $r\to R$.

Suppose we are interested in the asymptotics as $r\to R$ of the random counting
function $n_F(r)$, which
counts the number of zeroes of $F$ in the disk $\{|z|\le r\}$.
To simplify the notation, assume that $a_0=1$.
Denote by
\[
N_F(r) = \int_0^r \frac{n_F(t)}{t}\, {\rm d}t
\]
the integrated counting function. Then, by Jensen's formula,
\[
N_F(r) = \int_\bT \log|F(rt)|\, {\rm d}m(t) - \log|F(0)|
= \log\sigma_F(r) + \int_\bT \log|\widehat{F}_r(t)|\, {\rm d}m(t)
- \log|\xi_0|
\]
where $ \widehat{F}_r(t) = F(rt)/\sigma_F(r) $. Note that
\[
\widehat{F}_r(e^{2\pi {\rm i}\theta}) = \sum_{k\ge 0}\, \xi_k
\widehat{a}_k(r) e^{2\pi {\rm i}k\theta}
\]
is a random Fourier series normalized by the condition
$\sum_{k\ge 0} | \widehat{a}_k(r) |^2 = 1$.

\medskip
First, assume that the
$\xi_k$'s are standard complex Gaussian random variables. Then, for every $t\in\bT$,
the random variable $ \widehat{F}_r(t) $ is again standard complex-valued and Gaussian, and
$ \cE \bigl| \log |\widehat{F}_r(t)| \bigr| $ is a positive numerical constant.
Therefore,
\[
\sup_{r<R}\, \cE \bigl| N_F(r) - \log\sigma_F(r) \bigr| \le C\,.
\]
Since both $ N_F(r) $ and $ \log \sigma_F(r) $ are convex functions, we can derive
from here that their derivatives are also close on average, i.e., that
\[
\cE \Bigl| n_F(r) - \frac{{\rm d}\log\sigma_F(r)}{{\rm d}\log r} \Bigr|
\]
is relatively small outside a small exceptional set $ E $ of values of $ r $ where
the derivative
\[
\frac{{\rm d} \log\sigma_F(r)}{{\rm d}\log r}
\]
changes too fast due to the irregular behaviour of $a_k$'s.
Invoking an appropriate version of the Borel-Cantelli lemma, we can also
establish an almost sure analogue of this result.

If we are interested in the
angular distribution of zeroes, the same idea works, we only need to replace Jensen's
formula by its modification for angular sectors.

\medskip
The same approach works for the Steinhaus coefficients
$\xi_k=e^{2\pi {\rm i} \gamma_k}$, where $\gamma_k$ are independent and
uniformly distributed on $[0, 1]$. In this case, one needs to
estimate the expectation of the modulus of the logarithm of the absolute value
of a normalized linear combination of independent Steinhaus
variables. This was done by Offord in~\cite{Offord-1968}; twenty
years later, Ullrich~\cite{Ullrich-88a, Ullrich-88b} and Favorov~\cite{Favorov-90, Favorov-94}
independently rediscovered this idea and gave new applications.

\medskip
A linear combination of Rademacher random variables $ x = \sum \xi_k a_k $ can vanish with
positive probability. Then one cannot hope to estimate from below the logarithmic expectation
$ \cE \bigl\{ \log |x | \bigr\} $.
In~\cite{L-O}, Littlewood
and Offord invented ingenious and formidable techniques to
circumvent this obstacle. Later, these techniques were further
developed by Offord in~\cite{Offord-1965, Offord-1995}, but unfortunately, they were
not sufficiently powerful to arrive to the same
conclusions as for the Gaussian and the Steinhaus coefficients.
Still, there is a reserve: note that in order to estimate the error
term in the Jensen formula we do not need to estimate
$\cE \bigl| \log|\widehat{F}_r(t)|\, \bigr| $ uniformly in $t\in\bT$.
For our purposes, the integral estimate of $\displaystyle \cE \Bigl\{
\int_\bT \bigl| \log|\widehat{F}_r(t)|\, \bigr|\, {\rm d}m(t)
\Bigr\}$ is not worse than the uniform bound for $\cE\bigl| \log|\widehat{F}_r(t)| \bigr|$.
To exploit this reserve,
we employ some harmonic analysis techniques.

\subsection{Logarithmic integrability of Rademacher Fourier series}
Let $\bigl( \Omega, \cP \bigr)$ be the probability space on which the Rademacher
random variables $\xi_k$ are defined. Denote by $Q=\Omega\times\bT$ the product
measure space with the product measure $\mu=\cP\times m$.

By $\RFfuncs \subset L^2(Q)$ we denote the closed subspace whose elements are
the Rademacher Fourier series (i.e., the closed linear span of
$\xi_k e^{2\pi {\rm i} k \theta}$), and $\|f\|_2$ always stands for the $L^2(Q)$-norm.

\medskip
Our first result is a distributional inequality, which says that if
a Rademacher Fourier series is small on a set $E\subset Q$ of
positive measure, then it must be small everywhere on $Q$.
\begin{theorem}\label{thm:Zygmund-type}
For each $f\in \RFfuncs$ and each set $E\subset Q$ of positive measure,
\[
\int_Q |f|^2\, {\rm d}\mu \le e^{C\log^6 (\frac2{\mu (E)})} \int_E
|f|^2\, {\rm d}\mu\,.
\]
\end{theorem}
The power $6$ on the RHS is not the best possible,
but we will show that it cannot be replaced by any number less than~$2$.
Note that this does not
contradict the possibility that the distributional
inequality can be improved if one is ready to discard an event of
small probability.

The proof of Theorem~\ref{thm:Zygmund-type} is based on ideas from harmonic analysis
developed by the first-named author in~\cite{Nazarov-AA,
Nazarov-preprint} to treat lacunary Fourier
series. It uses a Tur\'an-type lemma from~\cite[Chapter~1]{Nazarov-AA},
and the technique of small shifts introduced
in~\cite[Chapter~3]{Nazarov-AA}.

\medskip
Theorem~\ref{thm:Zygmund-type} immediately yields
the following $L^p(\mu)$-bound for the logarithm of the Rademacher Fourier series.
\begin{corollary}\label{cor:log-integr}
For each $f\in \RFfuncs$ with $\| f \|_2 = 1$, and
for each $p\ge 1$,
\[
\int_Q \bigl| \log |f| \bigr|^p\, {\rm d}\mu \le (Cp)^{6 p}\,.
\]
\end{corollary}
\noindent Note that even the case $p=1$ of this corollary is already non-trivial and new.

\subsection{The range of random Taylor series in the unit disk}
One of the consequences of the logarithmic integrability is the
answer to an old question from Kahane's
book~\cite[p.xii]{Kahane}:

\smallskip\par\noindent{\em Suppose that
\[
F(z) = \sum_{k\ge 0} \xi_k a_k z^k
\]
is a Rademacher Taylor series with the radius of
convergence $1$ and with
\[
\sum_{k\ge 0} |a_k|^2 = +\infty\,.
\]
Does the range $F(\mathbb D)$ fill the complex plane almost surely}?

\smallskip\par\noindent
We will prove this, and even more.
\begin{theorem}\label{thm:range}
Suppose
\[
F(z) = \sum_{k\ge 0} \zeta_k z^k\,,
\]
where $\bigl\{ \zeta_k \bigr\}_{k\ge 0}$
is a sequence
of independent complex-valued symmetric random
variables  satisfying the conditions
\[
\limsup_{k\to\infty} |\zeta_k|^{1/k} = 1
\qquad and \qquad
\sum_{k\ge 0} |\zeta_k|^2 = +\infty\, \qquad {\rm a.s.}\,.
\]
Then, {\rm a.s.},
\[
\sum_{w\colon F(w)=b} (1-|w|) = \infty \qquad \forall\, b\in\mathbb C\,.
\]
\end{theorem}
Note that if the series $\sum_{k\ge 0} |\zeta_k|^2$ converges, then
the function $F$ belongs to the Hardy space $H^2$, and therefore
its $b$-points obey the Blaschke condition
\[
\sum_{w\colon F(w)=b} (1-|w|) < \infty\,.
\]

Theorem~\ref{thm:range} has some history. In 1972,
Offord~\cite{Offord-1972} proved this result in the case when
$\zeta_k$ are uniformly distributed on the unit circle. The proof he
gave also works for the Taylor series with Gaussian coefficients; see
also Kahane~\cite[Section~12.3]{Kahane}. According to the
``reduction principle''~\cite[Section~1.7]{Kahane}, the special case
$\zeta_k=\xi_k a_k$, where $\xi_k$ are independent Rademacher random
variables and $a_k$ is a non-random sequence of complex numbers
such that $\limsup_k |a_k|^{1/k} = 1$ and $\sum_k |a_k|^2 = \infty$,
should yield the general case. In the Rade\-macher case, the result
was known under some additional restrictions on the growth of
the deterministic coefficients $a_k$. In 1981, Murai~\cite{Murai} proved
it assuming that $\liminf |a_k|>0$. Soon afterwards, Jacob and
Offord~\cite{Jacob-Offord} weakened this assumption to
\[
\liminf_{N\to\infty} \frac1{\log N} \sum_{k=0}^N |a_k|^2  > 0\,.
\]
To the best of our knowledge, since then there was no improvement.

\medskip
Curiously enough, even in the case when $\zeta_k = \xi_k a_k$ with the standard complex
Gaussian $\xi_k$'s, the question when $F(\mathbb D)=\mathbb C$ almost surely is not completely
settled. Recall that in~\cite{Murai-lacunar} Murai proved Paley's conjecture,
which states that if $F$ is a (non-random) Taylor series with Hadamard gaps
and with the radius of convergence $1$, then $F$ assumes every complex value
infinitely often, provided that $ \sum_{k\ge 0} |a_k| = +\infty$.
Therefore, the same holds for random Taylor series with Hadamard gaps.
Even the case of sequences $a_k$ with a regular behaviour remains open:
\begin{question}\label{quest:gaussian}
Suppose that the non-random sequence $\{a_k\}$ is decaying regularly and satisfies
\begin{equation}\label{eq:Kahane_reg}
\sum_{k\ge 0} |a_k|^2<\infty\,, \qquad \sum_{k\ge 0} \Bigl|
\frac{a_k}{\sqrt{k}}\Bigr| = \infty\,,
\end{equation}
and suppose that $\xi_k$ are independent standard Gaussian
complex-valued random variables. Does the range of the random Taylor
series $F(z)=\sum_{k\ge 0} \xi_k a_k z^k$ fill the
whole complex plane $\mathbb C$ {\rm a.s.}?
\end{question}

Note that convergence of the first series in \eqref{eq:Kahane_reg} yields
that, a.s., the function $F$ belongs to all
Hardy spaces $H^p$ with $p<\infty$. Moreover, by the Paley-Zygmund
theorem~\cite[Chapter~5]{Kahane}, a.s.,
we have $e^{\lambda |\widehat{F}|^2}\in L^1(\mathbb T)$  for every positive
$ \lambda $, where
$\widehat F$ denotes the non-tangential boundary values of $F$ on
$\mathbb T$. On the other hand, by Fernique's
theorem~\cite[Chapter~15]{Kahane}, divergence of the second series
in \eqref{eq:Kahane_reg} yields that, a.s., $F$ is unbounded in $\mathbb D$.

\bigskip
\centerline{* \quad * \quad *}
\medskip
It is worth mentioning that our techniques can be applied to some
other questions about the distribution of zeroes of random Taylor series
including the one about the angular distribution of zeroes of random
entire functions in large disks. We plan to return to that question
in a separate paper.

\section{Proof of the distributional inequality for Rademacher Fourier series}

\subsection{List of notation}\label{subsect:notation}\mbox{}

\smallskip
$\bT = \bR/\bZ$; we also identify $\bT$ with the interval $[0, 1)\subset\bR$;

\smallskip
$m$ either the Lebesgue measure on $\bT$ normalized by $m(\bT)=1$, or the Lebesgue measure on $\bR$;

\smallskip
$e(\theta) = e^{2\pi {\rm i} \theta}$,  $ \theta\in\bT $;

\smallskip
$\bR_+ = (0, \infty)$;

\smallskip
$(\Omega, \cP)$ a probability space;

\smallskip
$\xi_k\colon \Omega \to \{\pm 1\}$, $k\in\bZ$, independent
Rademacher random variables;

\smallskip
$(Q, \mu) = (\Omega\times\bT, \cP\times m)$ product measure space,
$L^2(Q)=L^2(Q, \mu)$;

\smallskip
$\phi_k (\omega, \theta) = \xi_k(\omega) e(k\theta)$, $k\in\bZ$,
$(\omega, \theta)\in Q$;

\smallskip $\RFfuncs \subset L^2(Q)$ the subspace of Rademacher Fourier series
$\displaystyle f\!=\!\sum_{k\in\bZ} a_k \phi_k$, $\displaystyle
\sum_{k\in\bZ}|a_k|^2\!<\!\infty$.

\medskip The system $\{ \phi_k \}$ is an orthonormal basis in the space $\RFfuncs$,
and for $f\in \RFfuncs$, we have
\begin{multline*}
\| f \|_2^2 = \int_Q |f|^2\, {\rm d}\mu = \int_\Omega |f(\omega,
\cdot)|^2\, {\rm d}\cP(\omega) \\
= \int_{\bT} |f(\cdot, \theta)|^2\,
{\rm d}m(\theta) = \sum_{k\in\bZ} |a_k|^2 = \| \{ a_k\}
\|^2_{\ell^2(\bZ)} \,.
\end{multline*}

\smallskip For a set $E\subset Q$, we denote its sections by
$ E_\omega \eqdef \bigl\{\theta\in\bT\colon (\omega, \theta)\in E \bigr\}$,
$ \omega\in\Omega $.

\smallskip The set $E\subset Q$ shifted by $t\in\bT$ is denoted by
 $E+t \eqdef \bigl\{(\omega, \theta)\colon (\omega, \theta-t)\in E \bigr\}$. Then
\[
E_\omega + t = \bigl\{\theta\colon \theta-t\in E_\omega \bigr\} =
(E+t)_\omega.
\]

\smallskip
We put $\Delta_t(E) \eqdef \mu \left( (E+t)\setminus E\right)$.

\smallskip The function $g\in L^2(Q)$ shifted by $t$ is denoted by $g_t$:
$g_t(\omega, \theta)= g(\omega, \theta+t)$.

\smallskip Note that for the indicator function of $E$, we have
$\bigl( \done_E \bigr)_t = \done_{E-t}$.

\smallskip A measurable function $ b $ on $Q$
that does not depend on $\theta$ will be called a {\em random constant}.

\smallskip We write $[x]$ for the integral part of $x$.

\subsection{The result}\label{subsect:result}
Here is the main result of this part of the paper. It shows that an
arbitrary function $f$ in $\RFfuncs$ cannot be too close to
a random constant $b$, provided that
the uniform norm of $b$ is small compared with the $L^2$-norm of
$f$. The version we gave in the introduction corresponds to the case
when $b$ is the zero function. The extension below is needed for the
proof of Theorem~\ref{thm:range} on the range of random Taylor series in the
unit disk.

\begin{theorem}\label{thm:Zygmund-type-final}
For each $f\in \RFfuncs$, for each  random constant $b\in L^\infty (\Omega)$ with
$\| b \|_\infty < \tfrac1{20} \| f \|_2$, and for each set $E\subset
Q$ of positive measure,
\[
\int_Q |f|^2\, {\rm d}\mu \le
\exp \left( C\log^6 \left(\frac2{\mu (E)}\right) \right) \int_E |f-b|^2\, {\rm d}\mu\,.
\]
\end{theorem}

As an immediate corollary, we get
\begin{corollary}\label{cor:log-integr-final}
For each $f\in \RFfuncs$ with $\| f \|_2 = 1$, for each $b\in
L^\infty (\Omega)$ with $\| b \|_\infty < \tfrac1{20}$, and for each
$p\ge 1$, we have
\[
\int_Q \bigl| \log |f-b| \bigr|^p\, {\rm d}\mu \le (Cp)^{6 p}\,.
\]
\end{corollary}

We note that the condition on the function $b$ is a technical one. Its purpose is to
avoid degenerate cases, for example, the case
when the functions $f$ and $b$ are both equal to $\xi_0$.

\subsection{The basic tools}\label{subsect:basic-tools}
Here is the list of the tools we will be using in the proof of
Theorem~\ref{thm:Zygmund-type-final}.

\subsubsection{}\label{subsubsect:Turan}{\bf Tur\'an-type lemma}
~\cite[Chapter~I]{Nazarov-AA}. {\em Suppose
\[
p(z) = \sum_{k=0}^n a_k e^{i \la_k t}, \qquad a_k\in\bC, \quad \la_0 <
\,...\, < \la_n \in \bR,
\]
is an exponential polynomial. Then for any interval $J\subset \bR$
and any measurable subset $E\subset J$ of positive measure,
\[
\max_J |p| \le  \Bigl( \frac{C m(J)}{m(E)} \Bigr)^{n} \sup_{E}
|p|\,.
\]
}

We will also use the $L^2$-bound that follows
from this estimate, see ~\cite[Chapter~III, Lemma~3.3]{Nazarov-AA}.
It states that under the same assumptions,

\begin{equation}
\| p \|_{L^2(J)} \le  \Bigl( \frac{C m(J)}{m(E)} \Bigr)^{n+\frac12}
\| p \|_{L^2(E)}\,.\label{eq:turan_lt_est}
\end{equation}

\subsubsection{}\label{subsubsect:Khinchin}{\bf Khinchin's inequality}.
{\em Let $\bigl\{\xi_k \bigr\}$ be independent Rademacher random variables,
and let $\bigl\{ a_k \bigr\}$
be complex numbers. Then for each $p\ge 2$, we have
\[
\Bigl( \cE \Bigl| \sum_k a_k \xi_k \Bigr|^p \Bigr)^{1/p}
\le C\sqrt{p} \Bigl( \sum_k |a_k|^2 \Bigr)^{1/2}.
\]
}

\subsubsection{}\label{subsubsect:bilinear-Khinchin}{\bf Bilinear Khinchin's inequality}.
{\em
Let $\bigl\{\xi_k \bigr\}$ be independent Rademacher random variables,
and let $\bigl\{ a_{k, \ell} \bigr\}$
be complex numbers. Then for each $p\ge 2$, we have
\[
\Bigl( \cE \Bigl| \sum_{k\ne \ell} a_{k, \ell} \xi_k \xi_\ell \Bigr|^p \Bigr)^{1/p}
\le C p \Bigl( \sum_{k\ne\ell} |a_{k, \ell}|^2 \Bigr)^{1/2}.
\]
} A simple and elegant proof of this inequality can be
found in a recent preprint by L.~Erd\H{o}s, A.~Knowles, H.-T.~Yau,
J.~Yin~\cite[Appendix~B]{EKYY}.

\subsection{The class $\rm Exp_{\tt loc}$ of functions with almost linearly
dependent small shifts}\label{subsect:Exp_loc} The proof of
Theorem~\ref{thm:Zygmund-type-final} uses the technique of small
shifts developed in~\cite[Chapter~III]{Nazarov-AA}. In this and
the next two sections we will outline this technique.

\medskip
Let $\cH$ be a Hilbert space.
By $L^2(\bT, \cH)$ we denote the Hilbert space
of square integrable $\cH$-valued functions on $\bT$ (in the sense of Bochner).
Note that the space $L^2(\bT, L^2(\Omega))$ can be identified with  $L^2(Q)$.
To define the class of functions in $L^2(\bT, \cH)$ with almost linearly
dependent shifts, we introduce the following set of parameters:

\smallskip\noindent$\bullet$ the order $n\in\bN$ (a large parameter);

\smallskip\noindent$\bullet$ the localization parameter $\tau>0$ (a small parameter);

\smallskip\noindent$\bullet$ the error $\varkappa>0$ (a small parameter).

\begin{definition}[$\rm Exp_{\tt loc}$]\label{def:Exp}
We say that a function $g\in L^2(\bT, \cH)$ belongs to the class
${\rm Exp}_{\tt loc} (n, \tau, \varkappa, \cH)$ if for each $t\in (0, \tau)$
there exist complex numbers $a_k = a_k(t), k \in \{ 0, \dots, n \}$, with
$\sum_{k=0}^n |a_k|^2=1$,
such that
\[
\Bigl\| \sum_{k=0}^n a_k g_{kt }\Bigr\|_{L^2(\bT, \cH)} < \varkappa\,.
\]
\end{definition}

In the case $\cH=\bC$, this class was introduced
in~\cite[Chapter~III]{Nazarov-AA}. ``In small'' (i.e., on intervals
of length comparable with $\tau$), the functions from this class
behave similarly to exponential sums with $n$ frequencies
and with coefficients in $\cH$.
On the other hand, since the translations act continuously in
$L^2(\bT, \cH)$, for any given $g\in L^2(\bT, \cH)$, $n\in\bN$, $\varkappa>0$,
one can choose the parameter $\tau>0$ so small that
$g\in {\rm Exp}_{\tt loc} (n, \tau, \varkappa, \cH)$.

In the next three sections, we
extend main results about this class (the spectral description, the local
approximability by exponential sums with $n$ terms, and the
spreading lemma) from the scalar case to the case considered here.
Since the proofs of these extensions are similar to the ones given
in~\cite{Nazarov-AA}, we relegate them to the appendices.

\subsection{Spectral description of the class $\rm Exp_{\tt loc}$}\label{subsect:approx-spectrum}
The first lemma shows that each function $g\in {\rm Exp}_{\tt loc} (n, \tau, \varkappa, \cH)$ has an
``approximate spectrum'' $\Lambda_g$, which
consists of $n$ frequencies so that
the Fourier transform of $g$ is small in the $\ell^2$-norm away from
these frequencies.

For $m\in\bZ$, $\Lambda\subset\bR$, let
\[
\theta_\tau (m)=\min (1, \tau |m|), \quad
\Theta_{\tau, \Lambda}(m) = \prod_{\lambda\in\Lambda} \theta_\tau (m-\la)\,.
\]
\begin{lemma}\label{lemma:approx-spectrum}
Given $g\in {\rm Exp}_{\tt loc} (n, \tau, \varkappa, \cH)$, there exists
a set $\Lambda = \Lambda_g \subset\bR$ of $n$ distinct frequencies such that
\[
\sum_{m\in\bZ} \| \,\widehat{g}(m)\, \|_{\cH}^2\, \Theta_{\tau,
\Lambda}^2(m) \le \bigl( Cn \bigr)^{4n} \varkappa^2\,.
\]
\end{lemma}

The proof of this lemma will be given in Appendix~A.

\subsection{Local approximation by exponential sums with $n$ terms}\label{subsect:local-approx}
Starting with this section, we assume that $\cH=L^2(\Omega)$. Then ${\rm Exp}_{\tt loc} (n, \tau, \varkappa,
L^2(\Omega) )\subset L^2(Q)$.

For a finite set $\Lambda\subset\bR$, denote by ${\rm Exp}(\Lambda, \Omega)$
the linear space of exponential polynomials with frequencies in
$\Lambda$ and with coefficients depending on $\omega$.
The next lemma shows that, for a.e. $\omega\in\Omega$, the function
$\theta\mapsto g(\omega, \theta)$, $g\in {\rm Exp}_{\tt loc} (n,
\tau, \varkappa, L^2(\Omega) )$, can be well approximated by
exponential polynomials from ${\rm Exp}(\Lambda, \Omega)$,
on intervals $J\subset [0, 1)$ of length comparable with $\tau$.

Suppose that $M > 1$ satisfies
\[ \ell = \frac{1}{\tau M}\in\bN\,,\]
and partition $\bT$ into $l$ intervals of length $M\tau$:
\[
\bT = \bigcup_{k=0}^{\ell-1} \Bigl[\frac{k}{\ell}, \frac{k+1}{\ell} \Bigr)\,.
\]

\begin{lemma}\label{lemma:loc-approx}
Let $M$ be as above and let $g\in {\rm Exp}_{\tt loc} (n, \tau, \varkappa, L^2(\Omega) )$.
There exists a non-negative function $\Phi\in L^2(Q)$ with
\[
\| \Phi \|_2 \le \bigl( Cn \bigr)^{2n} \varkappa,
\]
and with the following property:
\smallskip\noindent
for every interval $J\subset \bT$ in the partition, there exists an exponential
polynomial $p^J\in {\rm Exp}(\Lambda_g, \Omega)$ such that, for a.e.
$\omega\in\Omega$ and a.e. $\theta\in J$,
\[
\bigl| g(\omega, \theta) - p^J (\omega, \theta) \bigr| \le M^n \, \Phi (\omega, \theta)\,.
\]
\end{lemma}

The proof of this lemma will be given in Appendix~B.

\subsection{Spreading Lemma}\label{subsect:spreading}
The next lemma is the crux of the proof
of Theorem~\ref{thm:Zygmund-type-final}.

Given a set $E\subset Q$ of positive measure, we put $ \Delta_t(E) =
\mu\bigl( (E+t)\setminus E \bigr) $.

\begin{lemma}\label{lemma:spreading}
Suppose $g\in {\rm Exp}_{\tt loc} (n, \tau, \varkappa, L^2(\Omega) )$
and $E\subset Q$ is a set of positive measure. There exists a set
$\widetilde{E} \supset E$ of measure $\mu(\widetilde{E}) \ge \mu(E) + \tfrac12
\Delta_{n \tau} (E)$ such that, for each $b\in L^2(\Omega)$,
\[
\int_{\widetilde{E}} |g-b|^2\, {\rm d}\mu \le \Bigl( \frac{Cn^3}{\Delta_{n
 \tau}^2 (E)} \Bigr)^{2n+1} \left( \int_E |g-b|^2\, {\rm d}\mu +
\varkappa^2 \right)\,.
\]
\end{lemma}

This lemma follows from the previous lemma combined with the
Tur\'an-type estimate \eqref{eq:turan_lt_est}.
The proof of Lemma~\ref{lemma:spreading} will be given in
Appendix~C.

\subsection{Starting the proof of Theorem~\ref{thm:Zygmund-type-final}.
Zygmund's premise and the operator $A_E$}\label{subsect:Zygmund-idea}

Suppose that
\[
f = \sum_{k\in\bZ} a_k \phi_k\,, \qquad \phi_k(\omega, \theta) = \xi_k(\omega)e(k\theta)\,,\quad
\{a_k\}\in\ell^2(\bZ)\,,
\]
and that $b\in L^2(\Omega)$. Let $E\subset Q$ be a
measurable set of positive measure. Then
\begin{align*}
\int_E |f-b|^2\, {\rm d}\mu &= \int_E
\Bigl[ \sum_{k, \ell} a_k \bar a_\ell\, \phi_k \bar\phi_\ell
- 2\,\re{(f \bar b)}  + |b|^2 \Bigr]\, {\rm d}\mu \\
&\ge \int_E \Bigl[ \sum_k |a_k|^2 |\phi_k|^2 \Bigr] \, {\rm d}\mu
+ \int_E \Bigl[
\sum_{k\ne\ell} a_k \bar a_\ell\, \phi_k \bar \phi_\ell \Bigr] \, {\rm d}\mu
- 2\, \re{\langle f, \done_E b \rangle} \\
&= \mu (E) \| f \|^2_2 + \langle A_E f, f \rangle - 2\, \re{\langle f, \done_E b \rangle},
\end{align*}
where $A_E$ is a bounded self-adjoint operator on $\RFfuncs$, whose matrix
$\left( A_E(k,\ell) \right)_{k,\ell \in \bZ}$ in the orthonormal basis $\{ \phi_k \}$ is given by
\[
A_E(k, \ell) = \begin{cases}
\langle \done_E, \phi_k \bar \phi_\ell \rangle \,, & k\ne\ell; \\
0\,, & k=\ell\,.
\end{cases}
\]
To estimate the Hilbert-Schmidt norm of $A_E$, we observe that
the functions $\{ \phi_k \bar \phi_\ell \}_{k\ne\ell}$  form an orthonormal system in $L^2(Q)$, and that
each function from this system is orthogonal to the function $\done$.
Then
\[
\sum_{k\ne\ell} |A_E(k, \ell)|^2 + \bigl| \langle \done_E, \done  \rangle \bigr|^2
\le \| \done_E \|_2^2 = \mu (E)\,,
\]
and therefore,
\[
\| A_E \|_{HS} = \Bigl( \sum _{k\ne\ell} |A_E(k, \ell)|^2 \Bigr)^{1/2} \le
\sqrt{\mu (E) - \mu(E)^2}\,.
\]
This estimate is useful for sets $E$ of large measure.

\subsection{The sets $E$ of large measure}\label{subsect:large-meas}
For each $\mu \in (0,1)$, let $D(\mu) \in (1, +\infty]$ be the smallest value such that the inequality
\[
\int_Q |f|^2\, {\rm d}\mu \le D(\mu) \int_E |f-b|^2\, {\rm d}\mu
\]
is satisfied for every $E\subset Q$ with $\mu (E)\ge \mu$,
for every $f\in \RFfuncs$, and for every random constant $b\in L^\infty (\Omega)$ with
$\| b \|_\infty < \tfrac1{20} \| f \|_2$.

Using the estimates from the previous
section, we get
\begin{align*}
\int_E |f-b|^2\, {\rm d}\mu
& \ge \bigl( \mu (E) - \| A \| \bigr) \| f \|_2^2 - 2\, \| \done_E b \|_2 \, \| f \|_2 \\
& \ge \bigl( \mu (E) - \sqrt{\mu (E) - \mu (E)^2} - \tfrac1{10} \bigr) \| f \|_2^2
\ge \tfrac12\, \| f \|_2^2\,,
\end{align*}
provided that $\mu(E)\ge \tfrac9{10}$. That is, $D(\mu)\le 2$ for
$\mu\ge \tfrac9{10}$.

\medskip In order to get an upper bound for $D(\mu)$ for smaller values of $\mu$, first of
all, we need to get a better bound for the Hilbert-Schmidt norm
of the operator $A_E$.

\subsection{A better bound for the Hilbert-Schmidt norm of $A_E$}
\label{subsect:better-bound-HSnorm}
Here, using the bilinear Khinchin inequalities~\ref{subsubsect:bilinear-Khinchin},
we show that {\em for each $p\ge 1$,}
\[
\| A_E \|_{HS} \le Cp \cdot \mu(E)^{1-\frac1{2p}}\,.
\]
For sets $E$ of small measure, this bound is better than the one we gave
in~\ref{subsect:Zygmund-idea}.

\medskip\noindent{\em Proof:} First, using duality and then H\"older's inequality,
we get
\begin{align*}
\Bigl( \sum _{k\ne\ell} |A_E(k, \ell)|^2 \Bigr)^{1/2}
&=
\sup \Bigl\{ \Bigl| \sum_{k\ne\ell} A_E(k, \ell) g_{k, \ell} \Bigr|\colon
\sum_{k\ne \ell} |g_{k,\ell}|^2 \le 1 \Bigr\} \\
&=
\sup \Bigl\{ \Bigl| \int_Q \done_E \bar g\, {\rm d}\mu \Bigr|\colon
g\in\operatorname{span}\bigl\{ \phi_k \bar\phi_\ell\bigr\}_{k\ne\ell}\,, \ \| g \|_2\le 1 \Bigr\}\\
&\le \mu(E)^{1-\frac1{2p}} \cdot
\sup\Bigl\{ \| g \|_{2p}\colon  g\in\operatorname{span}\bigl\{ \phi_k \bar\phi_\ell\bigr\}_{k\ne\ell}\,,
\ \| g \|_2\le 1 \Bigr\}\,.
\end{align*}
Now, using the bilinear Khinchin inequality, we will bound $ \| g \|_{2p} $ by $C p \| g \|_2$.
Since $ g\in\operatorname{span}\bigl\{ \phi_k \bar\phi_\ell\bigr\}_{k\ne\ell} $,
\[
g(\omega, \theta) = \sum_{k\ne\ell} g_{k, \ell} \xi_k (\omega) \xi_\ell (\omega) e((k-l)\theta)\,,
\]
whence,
\begin{align*}
\int_Q |g|^{2p}\, {\rm d}\mu &= \int_{\bT} {\rm d} m(\theta)\,
\int_\Omega {\rm d}\cP(\omega)
\Bigl|  \sum_{k\ne\ell} g_{k, \ell} \xi_k (\omega) \xi_\ell (\omega) e((k-l)\theta) \Bigr|^{2p} \\
&\le \int_{\bT} {\rm d} m(\theta)\, (Cp)^{2p}
\Bigl( \sum_{k\ne \ell} \bigl| g_{k, \ell} e((k-\ell)\theta) \bigr|^2 \Bigr)^{p} \\
&= (Cp)^{2p} \Bigl( \sum_{k\ne \ell} |g_{k, \ell}|^2 \Bigr)^{p} = (Cp)^{2p} \| g \|_2^{2p}\,,
\end{align*}
completing the proof. \hfill $\Box$

\subsection{The subspace $V_{E, b}$}\label{subsect:V_E}
Let $p \ge 1$. We now show that there exists a positive numerical constant $C'$
with the following property.
{\em If $E\subset Q$ is a set of positive measure
and $b\in L^2(Q)$, then there exists a subspace $V_{E, b}\subset
\RFfuncs$ of dimension at most
\[
n = \Bigl[\, \frac{C'p^2}{\mu(E)^{1/p}} \,\Bigr]
\]
such that for each function $g \in \RFfuncs \ominus V_{E, b}$
and each $b_1=c\cdot b$ with $c\in\bC$, we have}
\[
\int_Q |g|^2\, {\rm d}\mu \le \frac2{\mu (E)}\, \int_E |g-b_1|^2\,
{\rm d}\mu\,.
\]

\medskip\noindent{\em Proof:} This result is a rather straightforward consequence of the estimates
from~\ref{subsect:Zygmund-idea} and~\ref{subsect:better-bound-HSnorm}.
We enumerate the eigenvalues of the
operator $A_E$ so that their absolute values form a non-increasing sequence: $|\sigma_1|\ge
|\sigma_2| \ge $ ... . Let $h_1$, $h_2$, ...  be the
corresponding eigenvectors. Let $m \in \bZ$ and denote by $\widetilde{V}_E$ the
linear span of $h_1$, ... , $h_m$. Then the norm of the restriction
$A_E $ to $L^2_{\tt RF}\ominus \widetilde{V}_E$ equals $|\sigma_{m+1}|$.
Therefore, if the function $g\in \RFfuncs \ominus \widetilde{V}_E$, then $
\bigl| \langle A_E g, g  \rangle \bigr| \le |\sigma_{m+1}|\cdot
\|g\|_2^2 $.

Next,
\begin{align*}
\sigma_{m+1}^2 &\le \frac1{m+1} \sum_{j=1}^{m+1} \sigma_j^2 \le
\frac1{m+1} \sum_{j=1}^\infty \sigma_j^2 \\
&= \frac1{m+1}\, \| A_E \|_{HS}^2 \le \frac{Cp^2}{m+1}\cdot \mu (E)^{2-\frac1p}
< \frac14\, \mu(E)^2\,,
\end{align*}
provided that
\[
m \ge \Bigl[\, \frac{C^\prime p^2}{\mu(E)^{1/p}} \,\Bigr] - 1
\]
and $C^\prime$ is chosen large enough.

Denote by $U_{E, b}$ the one-dimensional space spanned by the
projection of the function $\done_E \cdot b$ to $\RFfuncs$, and
put $V_{E, b}=\widetilde{V}_E + U_{E, b}$. Then, assuming that
$g\in \RFfuncs \ominus V_{E, b} \subset \RFfuncs \ominus \widetilde{V}_E$ and applying the estimate
from~\ref{subsect:Zygmund-idea}, we get
\begin{align*}
\int_E |g-b_1|^2\, {\rm d}\mu &\ge \mu (E) \| g \|^2_2 + \langle A_E g,
g \rangle - 2\, \re{\langle g, \done_E b_1 \rangle} \\
&\ge \mu (E) \| g \|_2^2 - \frac12  \mu (E) \| g \|_2^2  =
\frac{\mu(E)}2\, \| g \|_2^2.
\end{align*}
Since $\dim V_{E, b}$ is at most $\Bigl[ C^\prime p^2 \mu(E)^{-1/p} \Bigr]$, the proof is complete.
\hfill $\Box$

\medskip Note that it suffices to take $C' = 4 C^2 + 1$, where $C$ is the constant that appears in the
bilinear Khinchin inequality~\ref{subsubsect:bilinear-Khinchin}, though this is
not essential for our purposes.

\subsection{Placing $f\in \RFfuncs$ in the class $\operatorname{Exp}_{\tt loc}$.
Condition $(C_\tau)$}
\label{subsect:putting-f-into-Exp_loc}
Introduce the function
\[ n(p, \mu) \eqdef \bigl[\, C''p^2 \cdot \mu^{-\frac1p} \,\bigr] \]
where $C''>C'$ is a sufficiently large numerical constant.
Fix $p \ge 1$ and let $E\subset Q$ be a given set  of positive measure.
Put $n=n\bigl( p, \tfrac12 \mu(E) \bigr)$ and choose the small
parameter $\tau$ so that, for every $t\in (0, \tau]$,
\[
\mu \Bigl(  \bigcap_{k=0}^n (E-kt) \Bigr) \ge \frac12 \mu (E)\,.
\eqno (C_\tau)
\]
This is possible since the function
$t \mapsto \mu \Bigl(  (E-t) \bigcap E \Bigr)$ is continuous and equals $\mu(E)$ at $0$.

\medskip
Now we prove that {\em given a set $E\subset Q$ of positive measure, $b\in L^2(Q)$, and
$p\ge 1$, each function $f\in \RFfuncs$ belongs to
the class $\operatorname{Exp}_{\tt loc}(n, \tau, \varkappa, L^2(\Omega))$ with
\[
n = n(p, \tfrac12 \mu (E))\,,
\qquad \varkappa^2 = \frac{4(n+1)}{\mu (E)}\, \int_{E} |f-b|^2\, {\rm d}\mu\,,
\]
and arbitrary $\tau$ satisfying condition $(C_\tau)$.}

\medskip\noindent{\em Proof:}
To shorten the notation, we put
\[
E' = E'_t = \bigcap_{k=0}^n (E-kt)\,.
\]
Then for every $k \in \{ 0,\ldots,n \} \,,$
\[
\int_{E'} |f_{kt}-b|^2\, {\rm d}\mu \le
\int_{E-kt} |f_{kt}-b|^2\, {\rm d}\mu
= \int_E |f-b|^2\, {\rm d}\mu,
\]
since $b$ depends only on $\omega$, and so, $b_{kt} = b$.

Given $t\in (0, \tau]$, we choose $a_0, \ldots, a_n \in \bC$
with $\sum_{k=0}^n |a_k|^2 = 1$ so that the function $g = \sum_{k=0}^n a_k f_{kt}$
belongs to the linear space $L^2_{\tt  RF}\ominus V_{E',\, b}$.
This is possible since
\[ \dim V_{E',\, b} \le n(p, \mu (E')) \le n(p, \tfrac12 \mu(E)) = n\,.\]
Since the function $g$ is orthogonal to the subspace $V_{E',\, b}$, we
can control its norm applying the estimate from~\ref{subsect:V_E}
with $b_1 = b \cdot \sum_k a_k$:
\begin{align*}
\int_Q |g|^2\, {\rm d}\mu & \le \frac2{\mu (E')} \int_{E'} |g-b_1|^2\, {\rm d}\mu \\
&\le \frac4{\mu (E)}\int_{E'} \Bigl| \sum_{k=0}^n a_k \bigl( f_{kt}-b \bigr) \Bigr|^2\, {\rm d}\mu \\
&\le \frac4{\mu (E)}\int_{E'} \sum_{k=0}^n \bigl| f_{kt}-b \bigr|^2\, {\rm d}\mu \le
\frac{4(n+1)}{\mu (E)} \, \int_E |f-b|^2\, {\rm d}\mu\,.
\end{align*}
That is,
\[
\Bigl\| \sum_{k=0}^n a_k f_{kt} \Bigr\|_2 \le \varkappa\,,
\]
and we are done. \hfill $\Box$

\subsection{Spreading the $L^2$-bound. Condition $(C_E)$}\label{subsect:spreading-C_E}
We apply the spreading Lemma~\ref{lemma:spreading} to the function
$f$ and the set $E$. It provides us with a set $\widetilde E
\supset E$, such that $\mu (\widetilde E) \ge \mu (E) + \tfrac12 \Delta_{n
\tau} (E)$ and
\begin{align*}
\int_{\widetilde E} |f-b|^2\, {\rm d}\mu &\le
\Bigl( \frac{Cn^3}{\Delta_{n \tau}^2 (E)} \Bigr)^{2n+1}\,
\left( \int_E |f-b|^2\, {\rm d}\mu + \varkappa^2 \right) \\
&\le \Bigl( \frac{Cn^3}{\Delta_{n \tau}^2 (E)} \Bigr)^{2n+1} \cdot
\frac{C(n+1)}{\mu (E)}\,\int_E |f-b|^2\, {\rm d}\mu\,,
\end{align*}
where $ n=n\bigl( p, \tfrac12 \mu(E) \bigr) \le 2C''p^2 \cdot
\mu(E)^{-\frac1p} $. There is not much value in this spreading until we learn
how to control the parameter $\Delta_{n \tau} (E)$ in terms of our
main parameters $\mu (E)$ and $p$. Clearly, the bigger
$\Delta_{n\tau} (E)$ is, the better is our spreading estimate. Recall that till
this moment, our only assumption on the value of $\tau$ has been
condition $(C_\tau)$ at the beginning of section~\ref{subsect:putting-f-into-Exp_loc}.

Now we will need the following condition on our set $E$:
\[
\max_t \Delta_t(E) \ge \frac{1}{2n} \mu(E).
\eqno (C_E)
\]
If condition $(C_E)$ holds then we can find $\tau > 0$ such that
$\Delta_{n\tau}(E) = \tfrac1{2n}\, \mu (E)$, while for all
$t\in (0, n\tau)$, $\Delta_t (E) < \tfrac1{2n}\, \mu (E)$.

\medskip Such $\tau$ will automatically satisfy condition
$(C_\tau)$ used in the derivation of the spreading estimate. Indeed,
\begin{align*}
\mu \Bigl( \bigcap_{k=0}^n (E-kt) \Bigr) &= \mu \Bigl( E\setminus
\bigcup_{k=1}^n \bigl( E\setminus (E-kt) \bigr) \Bigr) \\
&\ge \mu (E) - \sum_{k=1}^n \mu \bigl( E \setminus (E-kt) \bigr) \\
&= \mu (E) - \sum_{k=1}^n \mu \bigl( (E+kt) \setminus E \bigr) \\
&= \mu (E) - \sum_{k=1}^n \Delta_{kt} (E) \ge \mu (E) - \sum_{k=1}^n
\frac{\mu (E)}{2n} = \frac12 \mu (E).
\end{align*}

It is easy to see that there are sets $E\subset Q$ of arbitrary small positive
measure that do not satisfy condition $(C_E)$.
We assume now that condition $(C_E)$ is satisfied, putting aside the
question ``{\em What to do with the sets $E$ for which $(C_E)$ does not\
hold}?'' till the next section.

Substituting the value $\Delta_{n\tau}=\tfrac1{2n}\, \mu(E)$ into the spreading
estimate and taking into account that $n\le 2C''p^2 \, \mu(E)^{-\frac1p}$, we finally get
\begin{align*}
\int_{\widetilde E} |f-b|^2\, {\rm d}\mu &\le \Bigl( \frac{Cn^5}{\mu
(E)^2} \Bigr)^{2n+1} \cdot \frac{n+1}{\mu (E)}\,\int_E |f-b|^2\,
{\rm d}\mu \\
&\le \Bigl( \frac{C p}{\mu (E)} \Bigr)^{Cp^2\mu (E)^{-1/p}} \,\int_E
|f-b|^2\, {\rm d}\mu\,,
\end{align*}
while
\[
\mu (\widetilde{E}) \ge \mu (E) + \frac{c}{p^2}\,\mu
(E)^{1+\frac1{p}}\,.
\]
This is the spreading estimate that we will use for the sets $E$
satisfying condition $(C_E)$.

\subsection{The case of sets $E$ that do not satisfy condition $(C_E)$}
\label{subsect:sets-without-C_E}
Now, let us assume that $E\subset Q$
is a set of positive measure that does not satisfy condition
$(C_E)$, that is, for each $t\in [0, 1]$,
$\Delta_t(E)<\tfrac1{2n}\, \mu (E)$. The simplest example is any set
of the form $E=\Omega_1 \times \bT$, $\Omega_1\subset \Omega$. For
these sets, $\Delta_t (E)=0$ for every $t$. We will show that this
example is typical, i.e., the sets $E$ that do not satisfy condition
$(C_E)$ must have sufficiently many ``long sections'' $E_\omega$.
More precisely, let
\[
\Omega_1 = \bigl\{\omega\in\Omega\colon m(E_\omega) > 1- \tfrac{1}{n} \bigr\} \,.
\]
We show that $ \cP \{\Omega_1\} > \tfrac12\, \mu (E) $.

\medskip\noindent{\em Proof}:
Let
\[
\Omega_2 = \Omega \setminus \Omega_1 =\bigl\{\omega\in\Omega\colon m(E_\omega) \le 1- \tfrac{1}{n} \bigr\}.
\]
Since condition $(C_E)$ is not satisfied, we have
\[
\int_0^1 \Delta_t(E)\, {\rm d}t < \frac1{2n}\, \mu (E)\,.
\]
A straightforward computation shows that
\[
\int_0^1 m\left(\left(E_\omega + t \right) \setminus E_\omega \right) \, {\rm d} t =
m(E_\omega) \left( 1 - m(E_\omega) \right).
\]
Since $m(E_\omega) \le 1- \tfrac{1}{n}$ implies that $m(E_\omega) \le n m(E_\omega) (1 - m(E_\omega))$, we get
\begin{align*}
\int_{\Omega_2} m(E_\omega) \, {\rm d}\cP(\omega)
&\le n \int_{\Omega_2} m(E_\omega) \left(1 - m(E_\omega)\right) \, {\rm d}\cP(\omega) \\
&\le n \int_{\Omega} m(E_\omega) \left(1 - m(E_\omega)\right) \, {\rm d}\cP(\omega) \\
&= n \int_0^1 \Delta_t(E) \, {\rm d} t < \frac12 \mu(E).
\end{align*}
Therefore,
\[
\cP \bigl\{ \Omega_1 \bigr\}  \ge \int_{\Omega_1} m(E_\omega) \, {\rm d}\cP(\omega) =
\mu(E) - \int_{\Omega_2} m(E_\omega) \, {\rm d}\cP(\omega) > \frac12 \mu(E).
\]
\hfill $\Box$

\noindent {\bf Remark:}
{\em Since} $n=n\left(p, \tfrac12 \mu(E) \right) \ge 2$ {\em if} $C^{\prime\prime} \ge 2$,
{\em we trivially have}
\[
\cP \bigl\{ \Omega_1 \bigr\}  \le \frac{n}{n-1} \int_{\Omega_1} m(E_\omega) \, {\rm d}\cP(\omega) \le
\frac{n}{n-1} \, \mu(E) \le 2 \, \mu(E).
\]

\subsection{Many ``long sections''}\label{subsect:many-long-sections}
Assume that the set $E$ does not satisfy condition $(C_E)$. We will show that
\[
\int_Q |f|^2 \, {\rm d} \mu \le \frac{4}{\mu (E)}\, \int_E |f-b|^2 \, {\rm d} \mu \,,
\]
where, as above, $b=b(\omega)$ is a random constant, $\| b \|_\infty < \tfrac1{20} \| f \|_2$.

Let $\mu = \mu (E)$ and $\Omega_1$ be as above.
We have
\begin{align*}
\int_E |f-b|^2 \, {\rm d} \mu & \ge \int_{\Omega_1}
\left( \int_{E_\omega} \, |f-b|^2  \, {\rm d} m \right) \, {\rm d}\cP(\omega) \\
& = \int_{\Omega_1}\, \int_\bT |f-b|^2 \, {\rm d} m \, {\rm d}\cP(\omega)
- \int_{\Omega_1}\,\int_{\bT \setminus E_\omega}\, |f-b|^2\, {\rm d} m\, {\rm d}\cP(\omega) = ({\rm I}) - ({\rm II})\,.
\end{align*}
Notice that by the result of section~\ref{subsect:spreading-C_E},
we have $2\mu\ge\cP\{\Omega_1\}\ge\tfrac{1}{2}\mu$.

Bounding integral $({\rm I})$ from below is straightforward: we have
\[
\int_\bT |f-b|^2 \, {\rm d} m \ge \bigl( \| f \|_2 - \| b \|_\infty \bigr)^2 \ge \frac{9}{10}\, \| f \|_2^2\,,
\]
whence,
\[
({\rm I}) \ge \frac9{10}\, \| f \|_2^2\, \cP \{ \Omega_1 \} \ge \frac{9}{20}\, \cdot\mu\, \| f \|_2^2\,.
\]

Now let us estimate the integral (II) from above. We have
\[
({\rm II}) \le 2  \int_{\Omega_1}\, \int_{\bT\setminus E_\omega} |f|^2 \, {\rm d} m \, {\rm d}\cP(\omega) +
2 \int_{\Omega_1} \int_{\bT\setminus E_\omega} |b|^2\, {\rm d} m\, {\rm d}\cP(\omega) = ({\rm II}_a) + ({\rm II}_b)\,.
\]
Estimating the second integral is also straightforward:
\[
({\rm II}_b) \le
\frac{4 \mu}{n} \, \| b \|_\infty^2
< \frac{1}{10}\, \mu\, \| f \|_2^2
\]
(recall that $n\ge 2$ and $\|b\|_\infty < \tfrac1{20}\|f\|_2$).
Furthermore,
\[
({\rm II}_a) = 2 \int_{\Omega_1}\, \int_\bT \done_{\bT\setminus E_\omega}\, |f|^2\, {\rm d} m\, {\rm d}\cP(\omega)
\le 2 \Bigl( \int_{\Omega_1}\, \int_\bT \done_{\bT\setminus E_\omega} \Bigr)^{\frac{1}{r}}
\Bigl( \int_{\Omega_1}\, \int_\bT |f|^{2s} \Bigr)^{\frac{1}{s}}
\]
with $\tfrac{1}{r} + \tfrac{1}{s} = 1$.
By Khinchin's inequality,
\[
\Bigl( \int_{\Omega}\, \int_\bT |f|^{2s} \Bigr)^{\frac{1}{s}} \le C s \, \|f\|_2^2\,.
\]
Hence,
\[
({\rm II}_a) \le \Bigl(\frac{2 \mu}{n}\Bigr)^{\frac{1}{r}}\, Cs \,\|f\|_2^2\,.
\]
Letting $\tfrac{1}{r} = \tfrac{p}{p+1}, \tfrac{1}{s} = \tfrac{1}{p+1}$ and recalling
that $n \ge \tfrac12 C''\,p^2 \mu^{-1/p}$ and that $p\ge 1$, we continue the estimate as
\[
({\rm II}_a) \le \Bigl( \frac{4 \mu^{1+\frac{1}{p}}}{C''\, p^2}\Bigr)^{\frac{p}{p+1}} \, 2 C p\, \|f\|_2^2
< \frac{8 C}{\sqrt{C''}}\, \mu \,  p^{-\frac{p-1}{p+1}} \, \| f \|_2^2 < \frac{1}{10}\, \mu\, \| f \|_2^2\,,
\]
provided that the constant $C''$ in the definition of $n$ was chosen sufficiently big.
Finally,
\[
\int_E |f-b|^2 \ge  ({\rm I}) - ({\rm II}_a) - ({\rm II}_b)
\ge \Bigl( \frac{9}{20} - \frac{4}{20}\Bigr) \mu \|f\|_2^2 = \frac{1}{4}\,\mu\, \| f \|_2^2\,,
\]
completing the argument. \hfill $\Box$

\subsection{End of the proof of Theorem~\ref{thm:Zygmund-type-final}: solving a difference
inequality}\label{subsect:difference_ineq}
Recall that by $D(\mu)$ we denote the smallest value such that the inequality
\[
\int_Q |f|^2\, {\rm d}\mu \le D(\mu) \int_E |f-b|^2\, {\rm d}\mu
\]
holds for every $E\subset Q$ with $\mu (E)\ge \mu$,
every $f\in \RFfuncs$, and every random constant $b\in L^\infty (\Omega)$ satisfying
$\| b \|_\infty < \tfrac1{20} \| f \|_2$.

By~\ref{subsect:large-meas}, $D(\mu) \le 2$ for $\mu\ge \tfrac9{10}$, and by
the estimates proven in~\ref{subsect:spreading-C_E} and~\ref{subsect:many-long-sections},
for $0<\mu<\tfrac9{10}$ we have
\[
D(\mu) < \max\Bigl\{ \Bigl( \frac{C p}{\mu} \Bigr)^{Cp^2\mu^{-\frac1{p}}}
D\bigl(\mu + \frac{c}{p^2}\,
\mu^{1+\frac1{p}}\bigr), \frac{4}{\mu} \Bigr\}\,.
\]
Increasing, if needed, the constant $C$ in the exponent,
and taking into account that $\frac{p}{\mu} \ge \frac{1}{\nicefrac{9}{10}} > 1$ and $D \ge 1$,
we simplify this to
\[
D(\mu) < \Bigl( \frac{p}{\mu} \Bigr)^{Cp^2\mu^{-\frac1{p}}}
D \left( \mu + \frac{c}{p^2}\, \mu^{1+\frac1{p}} \right).
\]
Put
\[
\delta (\mu) = \frac{c}{p^2}\, \mu^{1+\frac1p}\,.
\]
Making the constant $c$ on the right-hand side small enough,
we assume that $\delta\bigl( \tfrac9{10}\bigr)<\tfrac1{10}$
(it suffices to take $c<\tfrac1{10}$). Then, for $0<\mu<\tfrac9{10}$,
\[
\log D(\mu) - \log D(\mu +\delta (\mu) ) < Cp^2 \mu^{-\frac1p} \log\Bigl( \frac{p}\mu \Bigr)
< C\,\delta(\mu)\, p^4 \mu^{-1-\frac2{p}} \log\Bigl( \frac{p}\mu \Bigr).
\]
To solve this difference inequality, we define the sequence
$\mu_0 = \mu$, $\mu_{k+1} = \mu_k+\delta(\mu_k)$, $k \ge 0$, and stop when
$\mu_{s-1} < \tfrac9{10} \le \mu_s$. Since we assumed that
$\delta\bigl( \tfrac9{10} \bigr) < \tfrac1{10}$, the terminal value $\mu_s$
will be strictly less than $1$. We get
\begin{align*}
\log D(\mu) &= \log D(\mu_s) + \sum_{k=0}^{s-1}
\bigl[ \log D(\mu_k) - \log D(\mu_{k+1}) \bigr] \\
&< 1 + Cp^4\, \sum_{k=0}^{s-1} \delta(\mu_k) \mu_k^{-1-\frac2p} \log\Bigl( \frac{p}\mu_k \Bigr) \\
&<  1 + Cp^4\, \log\Bigl( \frac{p}\mu \Bigr)\,
\sum_{k=0}^{s-1} \delta(\mu_k) \mu_k^{-1-\frac2p}\,.
\end{align*}
Since
$ \mu_{k+1} = \mu_k + c\,p^{-2}\, \mu_k^{1+\frac1{p}} < C\mu_k $,
we have
$ \mu_k^{-1-\frac2{p}} < C \mu_{k+1}^{-1-\frac2{p}} $.
Therefore,
\begin{align*}
\sum_{k=0}^{s-1} \delta(\mu_k) \mu_k^{-1-\frac2p}
&< C\, \sum_{k=0}^{s-1} \delta(\mu_k) \mu_{k+1}^{-1-\frac2p} \\
&< C\, \sum_{k=0}^{s-1}
\int_{\mu_k}^{\mu_{k+1}} \frac{{\rm d}x}{x^{1+\frac2p}}
< C\, \int_{\mu}^1 \frac{{\rm d}x}{x^{1+\frac2p}}
< C\, p\, \mu^{-\frac2p}\,,
\end{align*}
whence,
\[ \log D(\mu) < 1 + C\, p^5\, \mu^{-\frac2p}\, \log\left( \frac{p}{\mu}\right). \]
This holds for any $p\ge 1$. Letting $p=2\log\bigl( \tfrac{2}{\mu} \bigr)$,
we finally get
$ \log D(\mu) < C \log^6\Bigl( \frac2\mu \Bigr) $.
This completes the proof of Theorem~\ref{thm:Zygmund-type-final}.
\hfill $\Box$

\section{Proof of Theorem~\ref{thm:range} on
the range of random Taylor series}\label{sect_Kahane}

First, we prove the theorem in the special case when $\zeta_n = \xi_n a_n$, where
$\xi_n$ are independent Rademacher random variables, and $\bigl\{ a_n \bigr\}$ is
a non-random sequence of complex numbers satisfying the conditions
$\limsup_n |a_n|^{1/n} = 1 $ and $\sum_n |a_n|^2 = \infty$. That part of the proof
is based on the logarithmic integrability of
the Rademacher Fourier series (Corollary~\ref{cor:log-integr-final}
to Theorem~\ref{thm:Zygmund-type-final}) combined with Jensen's formula.
Then using ``the principle of reduction'' as stated in the Kahane
book~\cite[Section~1.7]{Kahane}, we get the result in the general case.

\medskip Let us introduce some notation.
For $b\in\bC$, $0<r<1$ we denote by $n_F(r, b)$ the number of solutions to the
equation $F(z)=b$ in the disk $r\bD$, the solutions being counted with their
multiplicities. In this section it will be convenient to set
\[
N_F(r, b) \eqdef \int_{1/2}^r \frac{n_F(t, b)}{t}\, {\rm d}t\,.
\]
By Jensen's formula
\begin{equation}\label{eq:Jensen}
N_F(r, b) = \int_\bT \log |F(r e(\theta)) - b|\, {\rm d}m(\theta) -
\int_\bT \log |F(\tfrac12 e(\theta)) - b|\, {\rm d}m(\theta)\,.
\end{equation}
We will prove that a.s. we have
\[
\lim_{r \to 1} N_F(r, b) = \infty \,, \qquad \forall b\in\bC\,,
\]
which is equivalent to Theorem~\ref{thm:range}.

\subsection{Proof of Theorem~\ref{thm:range} in the Rademacher case}

We define the functions $\sigma_F$ and $\widehat{F}$ by
\[
\sigma^2_F(r) \eqdef \sum_{n\ge 0} |a_n|^2 r^{2n}, \qquad
\widehat{F}(z) \eqdef \frac{F(z)}{\sigma_F(|z|)}\,,
\]
and note that $\| \widehat{F}(r e(\theta)) \|_{L^2(\bT)} = 1$.

Let $M \in \bN$. For every $r \in (\tfrac12, 1)$,
the function $(\omega, b) \mapsto N_F(r, b)$ on
$\Omega \times \bC $ is measurable
in $\omega$ for fixed $b$ and continuous in $b$ for fixed $\omega$.
Therefore, we can find a measurable
function $b^* = b^*(\omega)$ such that $|b^*| \le M$ and
\[
\inf_{|b|\le M} N_F (r, b) \ge N_F(r, b^*) - 1.
\]
Then
\begin{multline*}
\inf_{|b|\le M} N_F(r, b) \ge \int_\bT \log |F(r e(\theta)) - b^*|\, {\rm d}m(\theta) -
\int_\bT \log |F(\tfrac12 e(\theta)) - b^*|\, {\rm d}m(\theta) - 1 \\
= ({\rm I}_1) - ({\rm I}_2) - 1\,.
\end{multline*}
Note that
\begin{equation*}
({\rm I}_2) \le \frac12 \log \left( \int_\bT |F(\tfrac12 e(\theta)) - b^*|^2\, {\rm d}m(\theta) \right)
\le \frac12 \log \left( 2 \sigma^2_F(\tfrac12) + 2 M^2 \right).
\end{equation*}
For the integral $({\rm I}_1)$, we have the following lower bound:
\begin{multline*}
({\rm I}_1)
= \log \sigma_F(r) + \int_\bT \log \bigl| \widehat F(r e(\theta)) - \sigma^{-1}_F(r)\cdot b^* \bigr|\, {\rm d}m(\theta) \\
\ge \log \sigma_F(r) - \int_\bT \Bigl| \log \bigl| \widehat F(r e(\theta)) - \sigma^{-1}_F(r)\cdot b^* \bigr|\,\Bigr|\,
{\rm d}m(\theta)\,.
\end{multline*}
If we assume that $r$ is so close to $1$ that $\sigma_F(r)\ge 20M$, then, using our result on the logarithmic integrability of the Rademacher Fourier series (Corollary~\ref{cor:log-integr-final}), we get
\begin{multline*}
\cP \Bigl\{ \int_\bT \Bigl| \log \bigl| \widehat F(r e(\theta)) - \sigma^{-1}_F(r)\cdot b^* \bigr|\,\Bigr|\, {\rm d}m(\theta) > T \Bigr\}  \\
\le \frac{1}{T}\, \cE \Bigl( \int_\bT \Bigl| \log \bigl| \widehat F(r e(\theta)) - \sigma^{-1}_F(r)\cdot b^* \bigr|\,\Bigr|\, {\rm d} m(\theta) \Bigr)
\le \frac{C}{T}\,,
\end{multline*}
for all $T>0$.

Taking $r=r_m$ so that $\log \sigma_F (r_m) = 2m^2$ and $T=m^2$,
and applying the Borel-Cantelli lemma, we see
that, for a.e. $\omega\in\Om$, there exists $m_0=m_0(\omega, M)$ such that,
for each $m\ge m_0$,
\[
\int_\bT \Bigl| \log \bigl|
\widehat F(r_m e(\theta)) - \sigma^{-1}_F (r_m)\cdot b^* \bigr|\,\Bigr|\, {\rm d}m(\theta) < m^2 \,,
\]
whence,
\[
\inf_{|b|\le M} N_F(r_m, b)
\ge m^2 - \frac12 \log \left( 2 \sigma^2_F(\tfrac12) + 2 M^2 \right) - 1, \qquad \forall m\ge m_0\,.
\]
Therefore, for every $M\in\bN$, there is a set $A_M\subset\Omega$ with $\cP(A_M)=1$ such that,
for every $\omega\in A_M$ and every $b\in\bC$ with $|b|\le M$, we have
\begin{equation}\label{eq:N(r,b)}
\lim_{r\to 1} N_F(r, b) = \infty\,.
\end{equation}
Let $ A=\bigcap_{M} A_M$. Then $\cP(A)=1$, and for
every $\omega\in A$, $b\in\bC$, we have~\eqref{eq:N(r,b)}.
Thus, the theorem is proved in the Rademacher case.

\subsection{Proof of Theorem~\ref{thm:range} in the general case}

For every $M\in\bN$, consider the event
\[
B_M = \Bigl\{ \omega\colon \lim_{r\to 1}\,\inf_{|b|\le M}\, N_F(r, b) = +\infty \Bigr\}.
\]
Given $r\in\bigl( \tfrac12, 1\bigr)$, the function
$\displaystyle \inf_{|b|\le M}\, N_F(r, b) $  is measurable in $\omega$
(note that the infimum here can be taken over any dense countable subset of
the disk $\{|b|\le M\}$).
Thus, the set $B_M$ is measurable and so is the set $B=\bigcap_M B_M$, and for every
$\omega\in B$, $b\in \bC$, we have~\eqref{eq:N(r,b)}.
It remains to show that $B$ holds almost surely.

To that end, we extend the probability space to $\Omega\times\Omega'$
and introduce a sequence of independent Rademacher random
variables $\bigl\{ \xi_n (\omega') \bigr\}$, $\omega'\in\Omega'$,
which are also independent from the random variables $\bigl\{ \zeta_n (\omega) \bigr\}$,
$\omega\in\Omega$, and consider the random analytic function
\[
G(z) = G(z; \om, \om') = \sum_{n\ge 0} \xi_n(\omega') \zeta_n(\omega) z^n\,,
\qquad (\om, \om')\in\Om\times\Om'\,.
\]
By the previous section, for fixed  $\zeta_n$'s (outside a set of
probability zero in $\Omega$), the event
\[
\left\{ \omega'\in\Omega'\colon \lim_{r\to 1}\, \inf_{|b|\le M} N_G(r, b) = +\infty \right\}
\]
occurs with probability $1$. Hence, by Fubini's theorem, the
event $B_M$ occurs a.s. and so does the event $B$.
Note that due to the symmetry of the distribution of $\zeta_n$'s,
the random variables $\{ \xi_n (\omega') \zeta_n(\omega) \}$
are equidistributed with $\{ \zeta_n(\omega) \}$. This yields the theorem in the general
case of symmetric random variables.
\hfill $\Box$
\label{sect:Kahane}

\section{An example}\label{sect_examples}

In this section, we will present an example that shows that the constant $6$
in the exponent on the RHS of the inequality proven in
Theorem~\ref{thm:Zygmund-type} cannot be replaced by any number smaller than $2$.

Let
\[
g_N(\theta) =
\left( \sin ( 2 \pi \theta ) \right)^{2N} = \left( \frac{e(\theta) - e(-\theta)}{2 i} \right)^{2N} =
\sum_{|n|\le N} a_n e(2 n \theta).
\]
The function $g_N$ satisfies
\begin{equation}\label{eq:g_N}
|g_N(\theta)| \le e^{-c N^2} \qquad \mbox{for} \,\, |\theta| \le e^{-CN},
\end{equation}
provided that $C$ is large enough.

Now consider the Rademacher trigonometric polynomial
\[
f_N (\theta) = \sum_{|n|\le N} \xi_n a_n e(2 n\theta)\,,
\]
denote by $X_N$ the event that $\xi_n=+1$ for all
$n\in \bigl\{ -N, ..., N\bigr\}$, and put $E_N=X_N \times T_N$, where
$T_N = [-e^{CN},e^{CN}] \subset \bT$ is the set from~\eqref{eq:g_N}.
Then \[ \mu (E_N) \ge 2^{-(2N+1)} \cdot e^{-CN} \ge e^{-CN}\,,\] while
\[
\int_{E_N} |f_N|^2\, {\rm d}\mu \le e^{-cN^2} \mu (E_N) \le e^{-cN^2}
\]
and
\[
\int_{\Omega\times\bT} |f_N|^2\, {\rm d}\mu = \int_\bT |g_N|^2 \, {\rm d} m \,.
\]
It is not difficult to see that the integral on the RHS is not less than $\tfrac{c}{N}$,
for some constant $c>0$. Recalling that $|\log \mu (E_N)|\le CN$, we see that
for every $\e>0$, $C>0$, the inequality
\[
\int_{Q} |f_N|^2\, {\rm d}\mu \le e^{C|\log \mu(E_N)|^{2-\e}} \int_{E_N} |f_N|^2\,
{\rm d}\mu
\]
fails when $N\ge N_0(\e, C)$.
This shows that one cannot replace $6$ by any number less than $2$. \hfill $\Box$
\label{sect:examples}



\appendix

\section{Proof of the approximate spectrum lemma~\ref{lemma:approx-spectrum}}

The proof of Lemma~\ref{lemma:approx-spectrum}, with small modifications, follows
\cite[Section~3.1]{Nazarov-AA}.
We start with the following observation: if $g\in L^2(\bT, \mathcal H)$ and $a_0(t)$, \ldots ,
$a_n(t)$ are complex numbers, then
the $m$-th Fourier coefficient of the function
\[
x \mapsto  \sum_{k=0}^{n}a_{k}(t)g_{kt}(x)
= \sum_{k=0}^{n}a_{k}(t)g(x+kt)\] equals
\[
\widehat{g}(m)\cdot\sum_{k=0}^{n}a_{k}(t)e(ktm)
=
\widehat{g}(m)\cdot q_{t}\left(e(tm)\right),
\]
where $q_{t}(z)={\displaystyle \sum_{k=0}^{n}a_{k}(t)z^{k}}$.
Slightly perturbing the coefficients $a_k(t)$, we may assume
without loss of generality that the coefficients $a_0(t)$
and $a_n(t)$ do not vanish for $0<t<\tau$ (so that, for every $t$ in this
range, the polynomial $q_t$ is exactly of degree $n$ and does not vanish at the origin)
and that the arguments of the roots of $q_t$ are all distinct.

By Parseval's theorem,
\begin{equation}\label{eq:g_t_dual}
\int_{\bT}\, \Bigl\| \sum_{k=0}^{n}a_{k}(t)g_{kt}(x) \Bigr\|_{\mathcal H}^2\,{\rm d} x
=
\sum_{m\in\bZ}
\bigl\| \widehat{g}(m) \bigr\|_{\mathcal H}^2\,
\bigl| q_t\left( e(tm) \right)\bigr|^{2}.
\end{equation}
If $g\in {\rm Exp}_{\tt loc} (n, \tau, \varkappa, \mathcal H)$, then we can choose $a_0, \ldots, a_k$ so that the LHS of
\eqref{eq:g_t_dual} will be small for each $t\in (0, \tau)$.
On the other hand, whenever the norm of $\widehat{g}(m)$ is large,
the RHS of \eqref{eq:g_t_dual} can be small only when $q_{t}(e(tm))$ is small.
The proof of Lemma~\ref{lemma:approx-spectrum} will be based on two facts.
The first is that, on average, $|q_t(e(tm))|$ is relatively large outside some exceptional set,
which can be covered by at most $n$ intervals of length $\tfrac{1}{4n(n+1) \tau}$.
The second is that there exists a $t_0$
such that $q_{t_0}(e(tm))$ can be effectively bounded from below on this exceptional set.

We start with a lemma on arithmetic progressions.
\begin{lemma}\label{lemma:arithm}
Given a measurable set $G\subset\bR_+$, put
\[
V_G = \Bigl\{ t\in \bigl(\tfrac12 \tau, \tau \bigr)\colon \exists k\in \bN
{\rm\ s.t.\ } \tfrac{k}{t}\in G \Bigr\}.
\]
Then $m(V_G) < \tau^2 m(G)$.
\end{lemma}
\noindent
This lemma shows that if $m(G)<\tfrac{1}{2 \tau}$, then there are significantly many points
$t\in (\tfrac12 \tau, \tau)$ such that no point $k/t$, $k\in\bN$, belongs to $G$.

\medskip\noindent{\em Proof of Lemma~\ref{lemma:arithm}}:
We have
\[
\sum_{k\in\bN} \done_G\bigl( \tfrac{k}{t} \bigr) \ge \done_{V_G} (t) \,.
\]
Integrating over $t\in\bigl( \tfrac12 \tau, \tau \bigr)$, we get
\begin{multline*}
m(V_G) \le \int_{\tau/2}^\tau \sum_{k\in\bN} \done_G\Bigl( \frac{k}{t} \Bigr)\,
{\rm d}t = \sum_{k\in\bN} k\, \int_{k/\tau}^{2k/\tau} \done_G(s)\, \frac{{\rm d}s}{s^2} \\
= \int_0^\infty \done_G(s) \Bigl( \sum_{s\tau/2<k<s\tau} k \Bigr) \frac{{\rm d}s}{s^2}
< \tau^2 \int_0^\infty \done_G(s)\, {\rm d}s = \tau^2 m(G)\,,
\end{multline*}
because $\sum_{s\tau/2<k<s\tau} k < \tau^2 s^2$.
\hfill $\Box$

\medskip The following lemma shows that the Fourier coefficients $\widehat{g}(m)$ are small
outside $n$ intervals of controlled length. Put
\[
\delta = \frac{1}{8 n(n+1)}\,.
\]
This choice of $\delta$ will stay fixed till the end of the proof of
Lemma~\ref{lemma:approx-spectrum}.
\begin{lemma}\label{lem:crit_intervals}
There exist $n$ intervals $I_{1},\ldots,I_{n}$
of length $\frac{2\delta}{\tau}$ such that
\[
\sum_{m\in\bZ\backslash\bigcup I_{j}}
\bigl\| \widehat{g}(m) \bigr\|_{\mathcal H}^2
< \Bigl( \frac{C}{\delta} \Bigr)^{2n}
\varkappa^{2}.
\]
\end{lemma}

\noindent{\em Proof of Lemma~\ref{lem:crit_intervals}}:
By the continuity of the shift in $L^2(\bT, \mathcal H)$,
we can assume that the coefficients $a_{k}\left(t\right)$ are piecewise constant
functions of $t$, and hence measurable. Then, we can integrate
Parseval's formula~\eqref{eq:g_t_dual} over the
interval $\left(0,\tau\right)$. Recalling that the LHS of~\eqref{eq:g_t_dual}
is less than $\varkappa^2$, we get
\[
\sum_{m\in\bZ} \bigl\| \widehat{g}(m) \bigr\|_{\mathcal H}^2\, \rho^2(m) < \varkappa^2\,,
\]
where
\[
\rho^2(m) = \frac1{\tau}\, \int_0^\tau |q_t(e(tm))|^2\, {\rm d}t\,.
\]
Introduce the set
\[
S = \left\{ m\in\bZ\colon \rho^2(m) < \frac1{4(n+1)} \Bigl( \frac{\delta}{A} \Bigr)^{2n} \right\}.
\]
Here and elsewhere in this section,
$A$ is the positive numerical constant from the RHS of the
Tur\'an-type Lemma~\ref{subsubsect:Turan}. Then
Lemma~\ref{lem:crit_intervals} will follow from the following claim:
\begin{multline}\label{eq:cond_S}
S {\rm \ cannot\ contain\ } n+1 {\rm\ integers\ } m_1 < \ldots < m_{n+1} \\
{\rm\ such\ that\ }
m_{j+1}-m_{j} > \frac{2\delta}{\tau} \,, \quad \forall j \in \{ 1, \ldots, n \}\,.
\end{multline}
Indeed, this condition yields that the set $S$ can be covered by at most $n$ intervals
$I_1$, \ldots, $I_n$ of length $2\delta/\tau$ and
\[
\rho^2 (m) \ge \frac1{4(n+1)} \Bigl( \frac{\delta}{A} \Bigr)^{2n}, \qquad
m\in\bZ\setminus \bigcup_j I_j\,,
\]
whence
\[
\sum_{m\in\bZ\backslash\bigcup I_{j}}
\bigl\| \widehat{g}(m) \bigr\|_{\mathcal H}^2
\le 4 (n+1) \Bigl( \frac{A}{\delta} \Bigr)^{2n} \varkappa^2
< \Bigl( \frac{C}{\delta}\Bigr)^{2n} \varkappa^2
\]
with some numerical constant $C$. Thus, we need to prove claim~\eqref{eq:cond_S}.

\medskip
Suppose that \eqref{eq:cond_S} does not hold, i.e., there are $n+1$ integers
$m_1 < \ldots < m_{n+1}$ with $m_{j+1}-m_j > 2\delta/\tau $ that belong to
the set $S$. Then
\begin{equation}\label{eq:aver}
\int_{\tau/2}^\tau \sum_{j=1}^{n+1} \bigl| q_t (e(tm_j)) \bigr|^2 \, {\rm d}t
< \frac{\tau}4 \Bigl( \frac{\delta}{A} \Bigr)^{2n}\,.
\end{equation}
We call the value $t\in \bigl( \tfrac12 \tau, \tau \bigr)$ {\em bad} if
\[
\sum_{j=1}^{n+1} \bigl| q_t (e(tm_j)) \bigr|^2 <
\Bigl( \frac{\delta}{A} \Bigr)^{2n}\,.
\]
Otherwise, the value $t$ is called {\em good}. By~\eqref{eq:aver},
the measure of good $t$'s is less
than $\tau/4$. In the rest of the proof we will show that the measure of bad $t$'s
is also less than $\tau/4$, and this will lead us to a contradiction, which
will prove Lemma~\ref{lem:crit_intervals}.

\medskip
We will use the following
\begin{claim}\label{claim:U}
Let  $q(z) = \sum_{k=0}^n a_k z^k$ with $\sum_{k=0}^n |a_k|^2 = 1$.
Given $\Delta \in (0,1)$, let
\[
U = \Bigl\{s\in \bT\colon \bigl| q(e(s)) \bigr| < \Bigl( \frac{\Delta}{A} \Bigr)^n \Bigr\}\,.
\]
Then the set $U$ is a union of at most $n$ intervals of length at most $\Delta$ each.
\end{claim}

\noindent{\em Proof of Claim~\ref{claim:U}}: $U$ is an open subset of $\bT$
which consists of open intervals (since $\Delta < 1$ and $A \ge 1$ we have that $U \ne \bT$ ).
The boundary points of these intervals satisfy the equation
$\bigl| q(e(s)) \bigr|^{2}=\left(\frac{\Delta}{A}\right)^{2n}$,
which can be rewritten as
\[
\left(\sum_{k=0}^{n}a_{k} z^{k}\right)
\left(\sum_{k=0}^{n}\overline{a_{k}} z^{-k}\right)
=\left(\frac{\Delta}{A}\right)^{2n}, \qquad
z=e(s)\,.
\]
The LHS of this equation is a rational function of degree at most $2n$, and therefore
the number of solutions is at most $2n$. Hence $U$ consists of $l \le n$ intervals
$J_1$, \ldots, $J_l$, $l\le n$.

Next, note that since the sum of squares of the absolute values of the coefficients of $q$ equals $1$,
we have $\displaystyle \max_{s\in\bT} |q(e(s)|\ge 1$. Then,
applying Lemma~\ref{subsubsect:Turan} to the exponential polynomial $s\mapsto q(e(s))$,
we get
\[
1 \le\sup_{s\in\bT}\bigl| q\left( e(s) \right) \bigr|
\le
\Bigl( \frac{A}{m(J_i)} \Bigr)^n \cdot \sup_{s\in J_i} \bigl| q(e(s)) \bigr|
\le \Bigl( \frac{\Delta}{m(J_i)} \Bigr)^n\,.
\]
Hence, $m(J_i)\le \Delta$, proving the claim. \hfill $\Box$

\medskip Note that in the proof of this claim we did not use the full strength of
Tur\'an's lemma. For instance, we could have used the much simpler Remez' inequality.

\medskip
Now for $t\in \bigl(\tfrac12 \tau, \tau \bigr)$ consider the set
\[
S_t =
\Bigl\{
m\in\bZ\colon
\left|q_{t}\left(e\left(tm\right)\right)\right|<\Bigl(\frac{\delta}{A}\Bigr)^{n}
\Bigr\}\,.
\]
By the previous claim (applied with $\Delta = \delta$), there are points $\xi_1, ..., \xi_n\in\bR$
(centers of the intervals $J_i$) such that,
for each $m\in S_t$, there exist $i \in \{ 1,\ldots, n \}$ and $l \in \bZ$ such that
\begin{equation}\label{eq:S_t}
|t m - l - \xi_i| < \tfrac12\, \delta\,.
\end{equation}
Suppose that the value $t$ is bad. Then the $n+1$ integers $m_1$, \ldots, $m_{n+1}$ belong
to the set $S_t$, and by the Dirichlet box principle, there are two of these integers,
say $m_{j'}$ and $m_{j''}$ with $j'<j''$, that satisfy~\eqref{eq:S_t} with the same
value $i$. Then for this pair $ |t (m_{j''}-m_{j'}) - k|<\delta $,
with some non-negative integer $k$. Thus,
\[
\Bigl| \frac{k}t - (m_{j''}-m_{j'}) \Bigr| < \frac{\delta}t < \frac{2\delta}\tau\,.
\]
Note that since $m_{j''}-m_{j'}>\tfrac{2\delta}\tau$, the integer $k$ must be positive.
We conclude that the set of bad values $t$ is contained in the set $V_G$, where
$G$ is the union of $\tfrac12 n(n+1)$ intervals of length $\tfrac{4 \delta}\tau$
centered at all possible differences $ m_{j''}-m_{j'}$ with $j''>j'$.
The measure of the set $G$ is $ \tfrac{n(n+1)}2 \cdot \tfrac{4\delta}\tau $, which,
due to the choice of $\delta$, equals $\tfrac1{4\tau}$. By Lemma~\ref{lemma:arithm},
$m(V_G) < \tau^2 m(G) \le \tfrac14 \tau$. Thus, the measure of the set
of bad $t$'s is also less than $\tfrac14 \tau$, which finishes off the proof of
Lemma~\ref{lem:crit_intervals}. \hfill $\Box$

\medskip\noindent{\em Proof of Lemma~\ref{lemma:approx-spectrum}}:
We need to find a set $\Lambda = \Lambda_g \subset \bR$ of $n$ frequencies such that
\[
\sum_{m\in\bZ} \| \,\widehat{g}(m)\, \|_{\mathcal H}^2\, \Theta_{\tau,
\Lambda}^2(m) \le \bigl( Cn \bigr)^{4n} \varkappa^2\,,
\]
where
\[
\Theta_{\tau, \Lambda}(m) = \prod_{\lambda\in\Lambda} \theta_\tau (m-\la)\,,
\qquad \theta_\tau (m)=\min (1, \tau |m|).
\]
By Lemma~\ref{lem:crit_intervals}, there exists a collection of $n$ intervals $\{ I_j \}$,
each of length $\tfrac{2\delta}\tau$, such that
\[
\sum_{m\in \bZ \setminus \bigcup_j I_j} \| \, \widehat{g}(m) \, \|_{\mathcal H}^2\, \Theta^2_{\tau,\Lambda} (m)
\overset{\Theta \le 1}{\le}
\sum_{m\in \bZ \setminus \bigcup_j I_j} \| \, \widehat{g}(m) \, \|_{\mathcal H}^2
\le \left( C n \right)^{4n} \varkappa^2\,.
\]
Therefore, it remains to estimate the sum
\[
\sum_{m\in \bigcup_j I_j} \| \, \widehat{g}(m) \, \|_{\mathcal H}^2\, \Theta^2_{\tau,\Lambda} (m)\,.
\]
By Parseval's identity~\eqref{eq:g_t_dual}, for every $t\in (0, \tau)$,
\[
\sum_{m\in\bigcup_j I_j} \| \,\widehat{g}(m)\, \|_{\mathcal H}^2\,
| q_t( e(tm) )|^2 < \varkappa^2\,.
\]
Hence, it suffices to show that {\em there exist a value} $t_0\in (0, \tau)$ {\em and
a set} $\Lambda$ {\em of} $n$ {\em real numbers such that}
$ | q_{t_0} (e(t_0 m)) | \ge \delta^n \Theta_{\tau, \Lambda} (m) $
{\em for every} $m\in\bigcup_j I_j$.

\medskip First, we bound the absolute value of the polynomial $q_t$ from below by
the absolute value of another polynomial $p$ whose zeroes are obtained from the zeroes of $q_t$ by the radial
projection to the unit circle.
\begin{claim}\label{clm:zero_normalization}
Let $z_{j}\ne0$ for $1\le j \le n$,
and let $g(z)={\displaystyle c \cdot\prod_{j=1}^{n}\left(z-z_{j}\right)}$
be a polynomial of degree $n$ such that
${\displaystyle \sup_{|z|=1}\left|g\left(z\right)\right|\ge1}$.
Let ${\displaystyle h\left(z\right)=\prod_{j=1}^{n}\left(z-\zeta_{j}\right)}$,
where $\zeta_j = z_j/|z_j|$.
Then, for every $z\in\bT$,
\[
\left|h\left(z\right)\right| \le 2^{n} \left|g\left(z\right)\right|\,.
\]
\end{claim}

\noindent{\em Proof of Claim~\ref{clm:zero_normalization}}:
The ratio $\left|\frac{z-\zeta_{j}}{z-z_{j}}\right|$ attains its
maximum on $\{ |z|=1 \} $ at the point $z=-\zeta_{j}$, where it is equal
to $\frac{2}{1+\left|z_{j}\right|}$. Therefore,
\[
\left|\frac{h\left(z\right)}{g\left(z\right)}\right|
\le
\frac{1}{\left|c\right|}\prod_{j=1}^{n}\frac{2}{1+\left|z_{j}\right|}.
\]
By our assumption, there is some $z^{\prime}$, $|z'|=1$, such that
$\left|g\left(z^{\prime}\right)\right|\ge1$. Hence,
\[
1 \le
\left| c \right|\prod_{j=1}^{n}\left|z^{\prime} + z_{j}\right|
\le\left| c \right|\prod_{j=1}^{n}\left(1+\left|z_{j}\right|\right).
\]
Overall, we have
\[
\left|h\left(z\right)\right| \le
2^{n}
\left|g\left(z\right)\right|\cdot\frac{1}{\left|c\right|}
\cdot\prod_{j=1}^{n}\frac{1}{1+\left|z_{j}\right|}
\le2^{n} \left|g\left(z\right)\right|,
\]
proving the claim. \hfill $\Box$

\medskip Recall that $\displaystyle \sup_{|z|=1} |q_t(z)| \ge 1$. Hence,
applying Claim~\ref{clm:zero_normalization}, we conclude that
$ |q_t(z)| \ge 2^{-n}\, |p_t(z)| $ for $ |z|=1 $,
where $p_t$ is a monic polynomial of degree $n$ with all its zeroes on
the unit circle.

\medskip To choose $t_0$, we consider $n$ intervals $\widetilde{I}_j$ of length
$4\delta \tau^{-1}$ with the same centers as the intervals $I_j$ of
Lemma~\ref{lem:crit_intervals}, and put $\widetilde{S} = \bigcup_j \widetilde{I}_j$.
Let $\widetilde{G} = \widetilde{S} - \widetilde{S}$ be the difference set, with
$m(\widetilde{G})\le 8\delta \tau^{-1} \cdot n^2$.
We call the value $t\in \bigl( \tfrac12 \tau, \tau \bigr)$ {\em bad} if
there exists an integer $k\ne 0$ such that $k/t \in \widetilde G$. Since the set
$\widetilde G$ is symmetric with respect to $0$, we can estimate the measure of bad
$t$'s by applying Lemma~\ref{lemma:arithm} to the set $\widetilde{G} \cap \bR_+$.
Then the measure of bad values of $t$ is less than
$\tau^2 \cdot \tfrac12\, m(\widetilde G) \le 4 \delta \tau \cdot n^2 < \tfrac12 \tau$,
since $\delta \cdot 8n^2 < 1$. Therefore, there exists at least one {\em good }
value $t_0\in \bigl( \tfrac12 \tau, \tau \bigr)$ for which {\em every arithmetic
progression with difference $t_0^{-1}$ has at most one point in $\widetilde S$}.
We fix this value $t_0$ till the end of the proof.

\medskip
To simplify notation, we put $p=p_{t_0}$. The zero set of the function $x \mapsto p(e(t_0 x))$
consists of $n$ arithmetic progressions
with difference $t_0^{-1}$. By the choice of $t_0$, at most $n$ zeroes of
this function belong to the set $\widetilde S$. We denote these
zeroes by $\la_1$, ..., $\la_\ell$, $\ell \le n$. If $\ell < n$, we choose $n-\ell$
zeroes $\la_{\ell+1}$, ..., $\la_n$ in $\bR\setminus\widetilde S$ so that
$\bigl\{ e(t_0\la_j) \bigr\}_{1\le j \le n} $ is a complete set of zeroes of the
algebraic polynomial $p$; we recall that these zeroes are all distinct.

It remains to define a set $\Lambda$ of $n$ numbers, and to estimate from
below $|p(e(t_0m))|$ when $m\in\bigcup_j I_j$. Denote by $d_j(m)$ the distance from the integer
$m$ to the nearest point in the arithmetic progression
$\bigl\{ \la_j + k t_0^{-1} \bigr\}_{k\in\bZ}$.
We have
\[
\bigl| p(e(t_0 m)) \bigr| = 2^n \prod^n_{j=1} \bigl| \sin ( \pi t_0 (m-\la_j) \bigr|
\ge 2^n \prod_{j=1}^n \bigl( 2 t_0\, d_j(m) \bigr) \ge
2^n \tau^n \prod_{j=1}^n d_j(m)\,.
\]
We put $\Lambda = \bigl\{ \la_j \bigr\}_{1\le j \le n}$. Recall that here $m\in\bigcup_j I_{j}$, $\widetilde S = \bigcup_j \widetilde I_j$,
and that the arithmetic progression
$\bigl\{ \la_j + k t_0^{-1} \bigr\}_{k\in\bZ}$ either misses the set $\widetilde S$,
or has at most one element in $\widetilde S$.
In the first case, we get $d_j(m) \ge \delta \tau^{-1}$, while in the second case,
$d_j(m) \ge \min \bigl\{ \tfrac{\delta}{\tau},\left|m-\lambda_{j}\right| \bigr\}$.
Therefore, in both cases,
\begin{equation*}
d_{j}\left(m\right)
\ge \min \Bigl\{ \frac{\delta}{\tau},\left|m-\lambda_{j}\right| \Bigr\}
\stackrel{\delta\le\frac12}\ge
\frac{\delta}{\tau}\min\bigl\{ 1,\tau\left|m-\la_{j}\right|\bigr\}
= \frac{\delta}{\tau}\cdot \theta_{\tau}(m-\la_{j})\,.
\end{equation*}

Tying the ends together, we get
\begin{multline*}
\bigl| q_{t_0} (e(t_0 m)) \bigr| \ge 2^{-n} \bigl| p(e(t_0 m)) \bigr|
\ge
2^{-n} \cdot 2^n \tau^n \prod_{j=1}^n d_j(m) \\
\ge \tau^n \cdot
\Bigl( \frac{\delta}\tau \Bigr)^n
\Theta_{\tau, \Lambda} (m) = \delta^n \Theta_{\tau, \Lambda} (m)\,.
\end{multline*}
This completes the proof of Lemma~\ref{lemma:approx-spectrum}. \hfill $\Box$

\section{Proof of the lemma~\ref{lemma:loc-approx} on the local approximation}

The proof of Lemma~\ref{lemma:loc-approx} is very close to the
proof of the corresponding result in~\cite[Section~3.2]{Nazarov-AA}.
We start with a lemma on solutions of ordinary differential equations
(cf. Lemma~3.2 in~\cite{Nazarov-AA}).
\begin{lemma}\label{lem:sol_of_diff_eq_in_L2_Q}
Let
\[
D = \prod_{j=1}^n e(\la_j x)\, \frac{\rm d}{{\rm d}x}\,
e(-\la_j x)\, \qquad \la_1, \ldots, \la_n \in\bR\,,\quad \la_i \ne \la_j \mbox{  for } i \ne j\,,
\]
be a differential operator of order $n\ge 1$, and let $J\subset [0, 1]$ be an interval.
Suppose that $f\in L^2(\Omega \times J)$ and, for a.e. $\omega\in\Omega$,
$x\mapsto f(\omega, x)$ is a $C^n(J)$-function satisfying the
differential equation $Df=h$ with $h\in L^2(\Omega \times J)$.

Then there exists an exponential polynomial $p$ with spectrum $\la_1$, \ldots, $\la_n$, such that, for a.e. $\omega\in\Omega$,
\[
\sup_{x\in J} |f(\omega, x) - p(\omega, x)|
\le m(J)^n\, \frac1{m(J)}\, \int_J |h(\omega, x)|\, {\rm d}x\,.
\]
\end{lemma}

\noindent{\em Proof of Lemma~\ref{lem:sol_of_diff_eq_in_L2_Q}}:
Let $\phi $ be a particular solution of the equation $D\phi = h$ constructed by repeated
integration:
\[
\phi = \Bigl( \prod_{j=1}^n  e(\la_j x)\, \mathcal J
e(-\la_j x)\Bigr) h
\]
where $\mathcal J$ is the integral operator
\[
\bigl( \mathcal J \psi \bigr) (\omega, x) = \int_a^x \psi (\omega, t)\, {\rm d}t
\]
and $a$ is the left end-point of the interval $J$. Then, for a.e. $\omega$,
\[
| \phi (\omega, x) | \le m(J)^n\, \frac1{m(J)}\int_J | h(\omega, x) |\, {\rm d}x\,.
\]
The function $f-\phi$ satisfies the homogeneous equation $D(f-\phi)=0$. Hence,
$ p = f - \phi $ is an exponential polynomial with coefficients depending on
$\omega$:
\[
p(\omega, x) = \sum_{j=1}^n c_j(\omega) e(\la_j x)\,.
\]
\hfill $\Box$

\bigskip
Now we turn to the proof of~Lemma~\ref{lemma:loc-approx}.
We fix a function $g\in {\rm Exp}_{\tt loc} (n, \tau, \varkappa,
L^2(\Omega) )$. By Lemma~\ref{lemma:approx-spectrum}, this function
has an ``approximate spectrum''
$\Lambda = \Lambda_g = \left\{ \la_{j}\right\} _{1\le j \le n}$ so that
\[
\sum_{m\in\bZ} \| \,\widehat{g}(m)\, \|_{L^2(\Omega)}^2\, \Theta_{\tau,
\Lambda}^2(m) \le \bigl( Cn \bigr)^{4n} \varkappa^2\,,
\]
with
\[
\Theta_{\tau, \Lambda}(m) = \prod_{\lambda\in\Lambda} \theta_\tau (m-\la)\,,
\qquad
\theta_\tau (m)=\min (1, \tau |m|)\,.
\]
We fix $M > 1$ so that $1/(M \tau)$ is a positive integer, and
partition $\bT$ into intervals $J$ of length $M \tau$.

Put
\[
I_{k}=\Bigl(\la_{k}-\tfrac{1}{\tau},\la_{k}+\tfrac{1}{\tau}\Bigr),
\
\widetilde{I}_k=\Bigl(\la_{k}-\tfrac{2}{\tau},\la_{k}+\tfrac{2}{\tau}\Bigl),
\
E_{0}=\bR\backslash\bigcup_{k=1}^{n}I_{k}\,,
\
E_{k}=I_{k}\backslash{\bigcup_{j=1}^{k-1}I_{j}}\,.
\]
The sets $E_k$, $0\le k\le n$, form a partition of the real line.
Accordingly, we decompose $g$ into the sum $g={\displaystyle \sum_{k=0}^{n}g_{k}}$,
where $g_{k}$ is the projection of $g$ onto the closed subspace of $L^2(Q)$
that consists of functions
with spectrum contained in $E_{k}$. For each $k = 0, \ldots, n$, we have
\begin{equation}\label{eq:g_k}
\sum_{m\in\bZ}
\left\| \widehat{g_k}(m) \right\|_{L^2(\Omega)}^{2}
\Theta_{\tau, \Lambda}^2(m) < ( Cn )^{4n} \varkappa^2 \eqdef \widetilde{\varkappa}^2 \,.
\end{equation}
Since, for $m\in E_0$, $\Theta_{\tau, \Lambda}^2(m) \equiv 1$,
we get $ \| g_{0} \|_{L^2(Q)}  \le\widetilde{\varkappa} $.

Now let $1\le k \le n $. Let $n_k$ denote
the number of points $\la_j$ lying in $\widetilde{I}_{k}$.
We define a differential operator $D_{k}$ of order $n_{k}$ by
\[
D_{k} \eqdef
\prod_{\la_{j}\in\widetilde{I}_k}
e(\la_j x)
\frac{\rm d}{{\rm d}x}\,
e(-\la_j x).
\]
The function $g_{k}(x)$ is a trigonometric polynomial with coefficients depending
on $\omega$, hence, for a.e. $\omega$, it is an infinitely differentiable function of $x$.
We set $h_{k} \eqdef D_{k}g_{k}$.
Note that this is a trigonometric polynomial with
the same frequencies as $g_{k}$:
\[
\widehat{h}_{k}(\omega, m) =
\left(2\pi {\rm i}\right)^{n_{k}} \widehat{g}_{k}(\omega, m)
\prod_{\la_{j}\in\widetilde{I}_k}
\left(m-\la_{j}\right)\,.
\]
Consequently,
\[
\bigl| \widehat{h}_{k}(\omega, m) \bigr|
= (2\pi)^{n_k}\, \bigl| \widehat{g}_{k}(\omega, m) \bigr|
\prod_{\la_{j}\in\widetilde{I}_k}
| m-\la_{j}|\,.
\]
In the product on the RHS, $m\in E_k\subset I_k$ and $\la_j\in\widetilde{I}_k$.
Recalling the definition of the function $\theta_\tau$, we see that
\[
|m-\la_j| \le \tfrac3{\tau}\, \theta_\tau (m-\la_j)
\qquad {\rm for\ } m\in I_k, \ \la_j\in\widetilde{I}_k\,.
\]
Therefore,
\[
\bigl| \widehat{h}_{k}(\omega, m) \bigr|
\le \Bigl( \frac{6\pi}\tau \Bigr)^{n_k} \bigl| \widehat{g}_{k}(\omega, m) \bigr|
\prod_{\la_{j}\in\widetilde{I}_k}  \theta_\tau (m-\la_j)\,.
\]
Note that for $m\in E_k$ and for $\la_j\in\bZ\setminus \widetilde{I}_k$, we have
$\theta_\tau (m-\la_j)=1$. Thus,
\[
\bigl| \widehat{h}_{k}(\omega, m) \bigr|
\le \Bigl( \frac{6\pi}\tau \Bigr)^{n_k}
\bigl| \widehat{g}_{k}(\omega, m) \bigr| \Theta_{\tau, \Lambda}(m)\,,
\qquad \omega\in\Omega\,,
\]
whence, recalling estimate~\eqref{eq:g_k}, we obtain
\[
\| h_k \|_{L^2(Q)} \le \Bigl( \frac{6\pi}\tau \Bigr)^{n_k}\, \widetilde{\varkappa}\,.
\]

Applying Lemma \ref{lem:sol_of_diff_eq_in_L2_Q} to an interval $J$ of length $M \tau$,
we obtain an exponential polynomial $p_{k}^{J}$ with spectrum consisting
of frequencies $\la_{j}\in\widetilde{I}_{k}$ and with coefficients depending on $\omega$,
such that, for every $x\in J$ and almost every $\omega\in\Omega$,
\[
\bigl| g_{k}(\omega, x)-p_{k}^{J}(\omega, x) \bigr|
\le
( M \tau )^{n_{k}}
\cdot\frac{1}{M \tau} \int_J |h_{k}(\omega, t)|\, {\rm d}t\,.
\]
We denote by
\[
\frM f(\omega, x) =
\sup_{L\colon x \in L} \frac{1}{m(L)}
\int_L \left|f(\omega, t)\right|\, {\rm d} t
\]
the Hardy-Littlewood maximal function. The supremum is taken over all intervals
$L\subset [0, 1]$ containing $x$, but it is easy to see that it is enough to
restrict ourselves to the intervals with rational endpoints, which allows us
to rewrite $\frM f$ as $\sup \bigl\{ F_{\alpha, \beta}\colon \alpha, \beta\in\bQ\bigr\}$,
where
\[
F_{\alpha, \beta}(\omega, x) = \done_{[\alpha, \beta]} (x) G_{\alpha, \beta}(\omega)
\quad {\rm and} \quad G_{\alpha, \beta}(\omega)
= \frac1{\beta-\alpha}\int_\alpha^\beta |f(t, \omega)|\, {\rm d}t\,.
\]
By the Fubini theorem, $G_{\alpha, \beta}$ are measurable functions on $\Omega$, so
$F_{\alpha, \beta}$ are measurable functions on $Q$ and, thereby,
$\frM$ is measurable on $Q$ as well.

Let $\widetilde{h}_{k}=\tau^{n_{k}}h_{k}$.
Then
\[
\left|g_{k}\left(\omega, x\right)-p_{k}^{J}\left(\omega, x\right)\right|
\le M^{n_{k}} \cdot \frM\widetilde{h}_k (\omega, x)
\stackrel{M > 1}\le
M^n \cdot \frM\widetilde{h}_k (\omega, x)\,.
\]
Using the classical estimate for the $L^2$-norm of the maximal
function, we get, for a.e. $\omega$,
\[
\int_{\bT}
\Big[ \frM\widetilde{h}_k (\omega, x) \Big]^2\, {\rm d}x
\le C\,
\int_{\bT} \Bigl| \widetilde{h}_k (\omega, x) \Bigl|^2\, {\rm d}x \,.
\]
Recalling that $\| \widetilde{h}_k \|_{L^2(Q)} < C^{n_k} \widetilde{\varkappa}$,
we obtain
\[
\bigl\| \frM\widetilde{h}_k \bigr\|_{L^2(Q)}^2 = \int_{\Omega\times\bT}
\Big[ \frM\widetilde{h}_k (\omega, x) \Big]^2\, {\rm d}x\,
{\rm d}\mathcal P(\omega)
\le C\, \| \widetilde{h}_k \|_{L^2(Q)}^2 \le C^{2n_k} \widetilde{\varkappa}^2\,.
\]

We now set $p^{J} \eqdef {\displaystyle \sum_{k=1}^{n}p_{k}^{J}}$.
Notice that all the frequencies of the polynomial $p^{J}$ belong to
the set $\Lambda_g$. Then, for every $x\in J$,
\begin{eqnarray*}
\left| g\left(\omega, x\right)-p^{J}\left(\omega, x\right) \right| &\le&
\left|g_{0}\left(\omega, x\right)\right| +
\sum_{k=1}^{n}\left|g_{k}\left(\omega, x\right)-p_{k}^{J}\left(\omega, x\right)\right| \\
&\le& \left|g_{0}\left(\omega, x\right)\right|
+ M^n \sum_{k=1}^{n} \frM\widetilde{h}_{k}(\omega, \theta)\\
&\le& M^n
\Bigl( \left|g_{0}\left(\omega, x\right)\right|
+ \sum_{k=1}^{n} \frM\widetilde{h}_{k}(\omega, x) \Bigr)
\eqdef M^n \Phi\left(\omega, x\right).
\end{eqnarray*}
It remains to bound the norm of the ``error function'' $\Phi$:
\[
\bigl\| \Phi \bigr\|_{L^2(Q)} \le \bigl\| g_0 \bigr\|_{L^2(Q)} +
\sum_{k=1}^n \bigl\| \frM\widetilde{h}_k \bigr\|_{L^2(Q)}
\le \widetilde{\varkappa} + \sum_{k=1}^n C^{n_k} \widetilde{\varkappa} \le C^n \widetilde{\varkappa}
\le \left( C n\right)^{2n}  \varkappa \,.
\]
This proves the desired result. \hfill $\Box$

\section{Proof of the spreading lemma}
Till the end of this section, we fix the function
$g\in {\rm Exp}_{\tt loc} (n, \tau, \varkappa, L^2(\Omega) )$,
the set $E\subset Q$ of positive measure, and the ``random constant''
$b\in L^2(\Omega)$.

We will use two parameters, $M > 1$, $\tfrac{1}{M\tau} \in \bN$ and $\gamma \in (0,1)$; their
specific values will be chosen later in the proof.

\medskip\noindent{\bf Definition:} Let $J$ be an interval of length $M \tau$ in the partition of $\bT$.
The interval $J$ is called {\em $\omega$-white} if $m(J\cap E_{\omega})\ge\gamma m(J)$;
otherwise it is called {\em $\omega$-black}.

\medskip
Given $\omega$, the union of all $\omega$-white intervals will be denoted by
$W_\omega$. By $W\subset Q$ we denote the union of all sets $W_\omega$. Similarly, we
denote by $B_\omega$ the union of all $\omega$-black intervals and by $B\subset Q$ the
union of all sets $B_\omega$. Since we can write the set $W$ as
\[
\bigcup_J \bigl\{\omega\colon m(J\cap E_\omega) \ge\gamma m(J) \bigr\} \times J
\]
and the function $\omega\mapsto m(J\cap E_\omega)$ is measurable on $\Omega$ for
every interval $J$ in the partition, we see that $W$ and $B=Q\setminus W$ are measurable subsets of
$Q$.

\medskip
Let $\Phi$ be the error function given by the Local Approximation Lemma.
The next lemma enables us to extend our estimates for $g-b$ from the set $E$ to the set $W$.
\begin{lemma}\label{lem:est_for_l2_f_in_white_parts}
We have
\[
\int_{W} \left|g-b\right|^{2}\, {\rm d}\mu
\le \left( \frac{C}{\gamma} \right)^{2n+1}
\left[
\int_{W\cap E} |g-b|^{2} \, {\rm d}\mu
+ M^{2n+1} \int_{W} \Phi^{2}\, {\rm d}\mu
\right].
\]
\end{lemma}

\noindent{\em Proof of Lemma~\ref{lem:est_for_l2_f_in_white_parts}}:
Let $J$ be one of the $\omega$-white intervals of length $M \tau$.
By~Lemma \ref{lemma:loc-approx}, for almost every $\omega\in\Omega$ and every
$\theta\in J$, we have
\[
\left|( g(\omega,\theta)-b(\omega) ) -
( p^{J}(\omega,\theta)-b(\omega))\right|
=\left|g\left(\omega,\theta\right)-p^{J}\left(\omega,\theta\right)\right|
\le M^n\Phi(\omega,\theta),
\]
where $p^{J}$ is a exponential polynomial with $n$ frequencies and coefficients
depending on $\omega$.
Therefore,
\begin{equation}
\int_{J}\left|g-b\right|^{2}\, {\rm d}{\theta}
\le 2\left(
\int_{J} \left|p^{J}-b\right|^{2}\, {\rm d}{\theta}
+ M^{2n} \int_{J} \Phi^{2}\, {\rm d}{\theta}
\right).
\label{eq:approx_for_f_min_b_on_J}
\end{equation}
Applying the $L^{2}$-version of the Tur\'an-type lemma to the exponential polynomial
$p^J-b$, which has at most $n+1$ frequencies, we get
\begin{align*}
\int_{J}
\left| p^{J} - b \right|^{2}\, {\rm d}{\theta}
&\le
\left( \frac{C\,m(J)}{m(J\cap E_\omega)}\right)^{2n+1}
\int_{J\cap E_\omega}
\left| p^{J} - b \right|^{2}\,
{\rm d}{\theta} \\
&\le
\left(\frac{C}{\gamma}\right)^{2n+1} \int_{J\cap E_\omega}
\left| p^{J} - b \right|^{2}\,
{\rm d}{\theta}\,.
\end{align*}
Plugging this into \eqref{eq:approx_for_f_min_b_on_J},
we find that
\[
\int_{J} \left|g-b\right|^{2}\, {\rm d}{\theta}
\le
\left(\frac{C}{\gamma}\right)^{2n+1}
\int_{J\cap E_\omega} \left|p^{J}-b\right|^{2}\, {\rm d}{\theta}
+ 2 M^{2n} \int_{J}\Phi^{2}\, {\rm d}{\theta}\,.
\]
Summing these estimates over all $\omega$-white intervals $J$, and using that
\[
| p^J - b | \le |g - b| + |g - p_J| \le |g - b| + M^n \Phi \,,
\]
we get
\[
\int_{W_\omega}
\left| g - b \right|^{2}\, {\rm d}{\theta}
\le
\left(\frac{C}{\gamma}\right)^{2n+1}
\left[
\int_{W_\omega \cap E_\omega}
\left| g - b \right|^{2}\, {\rm d}{\theta}
+ M^{2n} \int_{W_\omega} \Phi^{2}\, {\rm d}{\theta}\right]
\]
Integrating over $\omega$, we get the result. \hfill $\Box$

\medskip The effectiveness of this lemma depends on the size of the set $W \cap E^{\rm c}$.
The following lemma is very similar to Lemma~3.4 from~\cite{Nazarov-AA}.
For the reader's convenience, we reproduce its proof. Recall that
$\Delta_{n\tau}(E) = \mu\left((E+n\tau) \setminus E \right)$.

\begin{lemma}\label{lem:white_part_low_bnd}
For $\gamma < \tfrac12$,
\[
\mu ( W \cap E^{\rm c} ) \ge \Delta_{n\tau}(E) - \left(\gamma + \frac{n}{M}\right).
\]
\end{lemma}

\noindent{\em Proof of Lemma~\ref{lem:white_part_low_bnd}}:
We have
\begin{eqnarray*}
m\bigl( (E_\omega+n\tau)\setminus E_\omega \bigr)
&=&
m\bigl( (E_\omega+n\tau) \cap E^{\rm c}_\omega \bigr) \\
&=&
m\bigl( (E_\omega+n\tau) \cap E_\omega^{\rm c} \cap W_\omega \bigr)
+ m\bigl( (E_\omega+n\tau) \cap E_\omega^{\rm c} \cap B_\omega \bigr) \\
&\le&
m\bigl( W_\omega \cap E_\omega^{\rm c} \bigr)
+ m\bigl( (E_\omega+n\tau) \cap E_\omega^{\rm c} \cap B_\omega \bigr).
\end{eqnarray*}
We need to estimate the second term on the RHS.

If the interval $J$ is $\omega$-black, then
\begin{multline*}
m\bigl( J \cap E_\omega^{\rm c} \cap (E_\omega+n\tau) \bigr)
\le m\bigl( J \cap (E_\omega+n\tau) \bigr)
\le m\bigl( J\setminus (J+n\tau) \bigr) + m\bigl( (E_\omega+n\tau) \cap (J+n\tau) \bigr) \\
\le n\tau + m\bigl( E_\omega \cap J\bigr)
< n\tau + \gamma m(J)
= \Bigl( \frac{n\tau}{m(J)} + \gamma \Bigr) m(J)\,.
\end{multline*}
Summing this inequality over all $\omega$-black intervals $J$, and recalling that
$m(J)=M \tau$, we obtain
\[
m\bigl( (E_\omega+n\tau) \cap E_\omega^{\rm c} \cap B_\omega \bigr)
\le \Bigl(  \frac{n \tau}{M \tau} + \gamma \Bigr) \cdot m(B_\omega) \le \frac{n}{M} + \gamma \,.
\]
Integrating over $\Omega$ we get the required result. \hfill $\Box$

\medskip

\noindent{\em Proof of Lemma~\ref{lemma:spreading}}:
We write $\Delta = \Delta_{n \tau} (E)$ and put
\[ M_1 = \frac{8 n}{\Delta}. \]
We consider two cases, according to whether $M_1 \tau \le 1$ or not.

In the first case, we choose $M \in \left[ M_1, 2 M_1 \right]$, so that $1/(M \tau)$ is an integer.
Notice that $M>1$. We set $\gamma = \tfrac18 \Delta < \tfrac12$ and let
$\widetilde{E} = E \cup \left(W \cap E^{\rm c} \right) = E \cup W$,
where $W$ is the union of the corresponding white intervals.
By Lemma~\ref{lem:white_part_low_bnd},
\begin{equation*}
\mu ( W \cap E^{\rm c} )
\ge \Delta - \Bigl( \gamma + \frac{n}{M} \Bigr)
\ge \Delta - \left(\frac{\Delta}{8} + \frac{\Delta}{8}\right) > \frac{\Delta}{2}\,.
\end{equation*}
Furthermore, using Lemma~\ref{lem:est_for_l2_f_in_white_parts},
we get
\[
\int_{W} \left|g-b\right|^{2}\, {\rm d}{\mu}
\le \left(\frac{C}{\gamma}\right)^{2n+1} \left[ \int_{W\cap E} \left|g-b\right|^{2}\, {\rm d}{\mu}
+ M^{2n} \int_{W} \,\Phi^{2}\, {\rm d}{\mu} \right].
\]
Plugging in the values of the parameters $\gamma$ and
$M$ and taking into account the bound on the norm of $\Phi$,
we find that the RHS is
\begin{eqnarray*}
&\le &\left(\frac{C}{\Delta}\right)^{2n+1} \left[ \int_{W \cap E} \left|g-b\right|^{2}\, {\rm d}{\mu} +
\left(\frac{C\,n}{\Delta}\right)^{2n}\, \int_{W}\Phi^{2}\, {\rm d}{\mu} \right] \\
&\le &\left(\frac{C}{\Delta}\right)^{2n+1} \left[ \int_{E} \left|g-b\right|^{2}\, {\rm d}{\mu} +
\left(\frac{C\,n^3}{\Delta}\right)^{2n}\, \varkappa^{2} \right] \\
&\le &\left(\frac{C\,n^3}{\Delta^2}\right)^{2n+1} \left[ \int_{E} \left|g-b\right|^{2}\, {\rm d}{\mu} + \,
\varkappa^{2} \right].
\end{eqnarray*}

Now we consider the second case, when $M_1 \tau > 1$. We set $M = \tfrac{1}{\tau}$
(that is, there is only one interval in the `partition') and note that
\[  M = \frac{1}{\tau} < M_1 = \frac{8n}{\Delta}\,. \]
We set $\gamma = \frac{\Delta}{2}$, and once again
$\widetilde{E} = E \cup \left(W \cap E^{\rm c} \right) = E \cup W$.
Similarly to the first case, Lemma~\ref{lem:est_for_l2_f_in_white_parts} gives us
\[
\int_{W} \left|g-b\right|^{2}\, {\rm d}{\mu}
\le \left(\frac{C\,n^3}{\Delta^2}\right)^{2n+1} \left[ \int_{E} \left|g-b\right|^{2}\, {\rm d}{\mu}
+ \, \varkappa^{2} \right].
\]
We now show that there are sufficiently many $\omega$-white intervals that contain a noticeable
portion of $E^{\rm c}$.
We define the function $\delta(\omega) = m\left((E_\omega + n \tau) \setminus E_\omega \right)$ and notice that
\[
\int_\Omega \delta(\omega) \, {\rm d}\mathcal P(\omega) = \Delta\,.
\]
Let $L = \bigl\{\omega\in\Omega\colon \delta(\omega) > \tfrac12 \Delta \bigr\}$. It is clear that
\[
\int_L \delta(\omega) \, {\rm d}\mathcal P(\omega) \ge \frac{\Delta}{2} \, .
\]
For $\omega \in L$ we have that $m(E^{\rm c}_\omega), m(E_\omega) \ge \delta(\omega) > \tfrac{\Delta}{2} = \gamma$,
and therefore $L \times \bT \subset W$. Thus $\left(L \times \bT \right) \cap E^{\rm c} \subset W \cap E^{\rm c}$ and
\[
m(W \cap E^{\rm c}) \ge m(\left(L \times \bT \right)
\cap E^{\rm c}) = \int_L m(E_\omega^{\rm c} ) \, {\rm d}\mathcal P(\omega)
\ge \int_L \delta(\omega) \, {\rm d}\mathcal P(\omega) \ge \frac{\Delta}{2} \,,
\]
proving the lemma. \hfill $\Box$


\end{document}